# GLOBAL $L_2$-SOLUTIONS OF STOCHASTIC NAVIER–STOKES EQUATIONS

By R. Mikulevicius and B. L. Rozovskii[1]

*University of Southern California, Los Angeles*

This paper concerns the Cauchy problem in $\mathbf{R}^d$ for the stochastic Navier–Stokes equation
$$\partial_t \mathbf{u} = \Delta \mathbf{u} - (\mathbf{u}, \nabla)\mathbf{u} - \nabla p + \mathbf{f}(\mathbf{u}) + [(\sigma, \nabla)\mathbf{u} - \nabla \tilde{p} + \mathbf{g}(\mathbf{u})] \circ \dot{W},$$
$$\mathbf{u}(0) = \mathbf{u}_0, \qquad \operatorname{div} \mathbf{u} = 0,$$
driven by white noise $\dot{W}$. Under minimal assumptions on regularity of the coefficients and random forces, the existence of a global weak (martingale) solution of the stochastic Navier–Stokes equation is proved. In the two-dimensional case, the existence and pathwise uniqueness of a global strong solution is shown. A Wiener chaos-based criterion for the existence and uniqueness of a strong global solution of the Navier–Stokes equations is established.

**1. Introduction.** In this paper we are concerned with the Cauchy problem for the stochastic Navier–Stokes equation

$$\begin{aligned}\partial_t \mathbf{u} &= \Delta \mathbf{u} - (\mathbf{u}, \nabla)\mathbf{u} - \nabla p + \mathbf{f}(\mathbf{u}) + [(\sigma, \nabla)\mathbf{u} - \nabla \tilde{p} + \mathbf{g}(\mathbf{u})] \circ \dot{W}, \\ \mathbf{u}(0) &= \mathbf{u}_0, \qquad \operatorname{div} \mathbf{u} = 0,\end{aligned} \qquad (1.1)$$

in $\mathbf{R}^d$ and some generalizations of this equation. In (1.1), $\dot{W}$ is a time derivative of a Hilbert space-valued Brownian motion (e.g., space-time white noise) and the stochastic integral is understood in the Stratonovich sense. Here and throughout the rest of the paper, vector fields on $\mathbf{R}^d$ are denoted by boldface letters. This convention also applies if the entries of the vector field are taking values in a Hilbert space.

Received April 2003; revised August 2003.
[1]Supported in part by NSF Grant DMS-98-02423, ONR Grant N00014-03-1-0227 and ARO Grant DAAD19-02-1-0374.

Equation (1.1) stems from the dynamics of fluid particles given by the stochastic flow map

$$(1.2) \qquad \dot{\eta}(t,x) = \mathbf{u}(t,\eta(t,x)) + \sigma(t,\eta(t,x)) \circ \dot{W}, \qquad \eta(0,x) = x,$$

with undetermined local characteristics $\mathbf{u}(t,x)$ and $\sigma(t,x)$. The generalized random field $\sigma(t,x) \circ \dot{W}$ models the turbulent part of the velocity field, while $\mathbf{u}(t,x)$ models its regular component. In [27] and [29] it was shown, following the classical scheme of the Newtonian fluid mechanics, that the regular component $\mathbf{u}(t,x)$ of the flow map satisfies (1.1).

Our interest in stochastic flows of the form (1.2) is related in part to the recent progress made on the turbulent transport problem (see, e.g., [14, 15, 20] and others). In these works, turbulent velocity field is modeled by a stationary isotropic Gaussian vector field $\mathbf{V}(t,x)$ with covariance

$$\mathbf{E}\mathbf{V}(t,x)\mathbf{V}(s,y) = \delta(t-s)C(x-y).$$

In the divergence-free case, the spatial covariance $C$ is defined by its Fourier transform

$$\widehat{C}(z) = \frac{C_0}{(d-1)(1+|z|^2)^{(d+\kappa)/2}}\left(I - \frac{zz^T}{|z|^2}\right),$$

where $C_0 > 0$ and $0 < \kappa < 2$.

The centered velocity field $\mathbf{V}(t,x) - \mathbf{E}\mathbf{V}(t,x)$ can be realized by way of its identification with the random vector field

$$(1.3) \qquad \sigma(x) \cdot \dot{W} = \sum_{k \geq 0} \sigma^k(x) \dot{w}_k(t),$$

where $\{\sigma^k, k \geq 1\}$ is an orthonormal basis in the reproducing kernel Hilbert space $H_C$ corresponding to the kernel function $C$ and $(w_k(t))_{k \geq 1}$ are independent one-dimensional Brownian motions. In fact, $H_C$ is the subset of divergence-free fields in the Sobolev space $H^{(d+\kappa)/2}(\mathbf{R}^d, \mathbf{R}^d)$. It can be shown (see, e.g., [21]) that each $\sigma^k$ is divergence-free and Hölder-continuous of order $\kappa/2$.

With this application in mind, we study (1.1) with nonsmooth coefficients and free forces. In particular, for the coefficient $\sigma$ it is sufficient to assume that it is bounded and divergence-free (in the sense of generalized functions).

The aim of this paper is to develop an analog of the Leray theory of $L_2$-solutions for the stochastic Navier–Stokes equation (1.1) and its generalizations.

To this end, in Section 2 we prove the existence of a global weak (i.e., martingale) solution of the Cauchy problem for the stochastic Navier–Stokes equation in $\mathbf{R}^d$, $d \geq 2$, such that

$$(1.4) \qquad \mathbf{E}\left[\sup_{s \leq T} |\mathbf{u}(s)|_2^2 + \int_0^T |\nabla \mathbf{u}(s)|_2^2 \, ds\right] < \infty.$$

In addition, we prove that if $d = 2$ and the free forces are Lipschitz-continuous with respect to $u$, there is a unique strong (pathwise) solution with property (1.4). In this case, we prove the convergence in probability of the approximating sequence.

Section 3 deals with Wiener chaos expansions for stochastic Navier–Stokes equations. In this section we derive a system of deterministic PDEs for projections of a solution of stochastic Navier–Stokes equations on the Hermite–Fourier basis in $L_2$-space of functions adapted to the filtration generated by $W$. This system is usually referred to as a propagator. We demonstrate that the existence (uniqueness) of a solution of the propagator is a necessary and sufficient (in certain sense) condition for the existence (uniqueness) of a strong (pathwise) solution of the related Navier–Stokes equation.

The existence and uniqueness of $L_2$-solutions of stochastic Navier–Stokes equations was studied by many authors (see, e.g., [2, 5, 6, 8, 9, 10, 11, 12, 25, 28, 33, 34, 35, 36]). There is also substantial literature on more regular solutions, invariant measures, Kolmogorov equations and other related topics that are beyond the scope of this paper.

The main novel elements in the present paper are as follows:

To the best of our knowledge, all previous results on martingale ($L_2$) solutions for equations similar to (1.1) were limited to bounded domains. An extension to an unbounded domain is not trivial since in the latter case the direct application of the compactness method, which is central to the proof, fails.

The existence of a global martingale solution of the Navier–Stokes equation with random forcing

$$(1.5) \qquad \partial_t \mathbf{u} = \Delta \mathbf{u} - (\mathbf{u}, \nabla)\mathbf{u} - \nabla p + \mathbf{f}(\mathbf{u}) + \mathbf{g}(\mathbf{u}) \circ \dot{W}$$

in unbounded domain was proved in ([7], Theorem 1.1). Equation (1.5) does not include a conceptually important term $(\sigma, \nabla)\mathbf{u} \circ \dot{W}$. Accordingly, it does not cover the case of turbulent flows (e.g., Kraichnan velocity) which is central to our paper. Also, the related results in [7] are limited to the case when $d = 2$ or 3 and all the moments of the initial condition are finite. It should be noted that [7] addresses a number of interesting issues (e.g., solutions in weighted spaces, equations driven by a homogeneous Wiener process) that are beyond the scope of the present paper.

As it was mentioned before, in contrast to previous work (see, e.g., [11, 25]), we do not assume any regularity of the coefficients.

In [26] and [29] it was shown (under more restrictive assumptions) that existence of a strong solution of a Navier–Stokes equation implies the existence of a solution of the propagator. However, the converse statement, which is in many ways more desirable, was not known previously even for linear equations.

Some results of the present paper were announced at the recent Trento meeting (see [28]).

We conclude this section with an outline of some notations that will be used in the paper.

Let us fix a separable Hilbert space $Y$. The scalar product of $x, y \in Y$ will be denoted by $x \cdot y$.

If $u$ is a function on $\mathbf{R}^d$, the following notational conventions will be used for its partial derivatives:

$$\partial_i u = \frac{\partial u}{\partial x_i}, \qquad \partial^2_{ij} = \frac{\partial^2 u}{\partial x_i \, \partial x_j}, \qquad \partial_t u = \frac{\partial u}{\partial t}$$

and

$$\nabla u = \partial u = (\partial_1 u, \ldots, \partial_d u)$$

and $\partial^2 u = (\partial^2_{ij} u)$ denotes the Hessian matrix of second derivatives. Let $\alpha = (\alpha_1, \ldots, \alpha_d)$ be a multi-index; then $\partial^\alpha_x = \prod_{i=1}^d \partial^{\alpha_i}_{x_i}$.

Let $C_0^\infty = C_0^\infty(\mathbf{R}^d)$ be the set of all infinitely differentiable functions on $\mathbf{R}^d$ with compact support.

For $s \in (-\infty, \infty)$, write $\Lambda^s = \Lambda^s_x = (1 - \sum_{i=1}^d \partial^2/\partial x_i^2)^{s/2}$.

For $p \in [1, \infty]$ and $s \in (-\infty, \infty)$, we define the space $H_p^s = H_p^s(\mathbf{R}^d)$ as the space of generalized functions $u$ with the finite norm

$$|u|_{s,p} = |\Lambda^s u|_p,$$

where $|\cdot|_p$ is the $L_p$ norm. Obviously, $H_p^0 = L_p$. Note that if $s \geq 0$ is an integer, the space $H_p^s$ coincides with the Sobolev space $W_p^s = W_p^s(\mathbf{R}^d)$.

If $p \in [1, \infty]$ and $s \in (-\infty, \infty)$, $H_p^s(Y) = H_p^s(\mathbf{R}^d, Y)$ denotes the space of $Y$-valued functions on $\mathbf{R}^d$ so that the norm $\|g\|_{s,p} = \||\Lambda^s g|_Y|_p < \infty$. We also write $L_p(Y) = L_p(\mathbf{R}^d, Y) = H_p^0(Y) = H_p^0(\mathbf{R}^d, Y)$. Let $C_0^\infty(Y)$ be the space of $Y$-valued infinitely differentiable functions on $\mathbf{R}^d$ with compact support.

Obviously, the spaces $C_0^\infty, C_0^\infty(Y), H_p^s(\mathbf{R}^d)$ and $H_p^s(\mathbf{R}^d, Y)$ can be extended to vector functions (denoted by boldface letters). For example, the space of all vector functions $\mathbf{u} = (u^1, \ldots, u^d)$ such that $\Lambda^s u^l \in L_p$, $l = 1, \ldots, d$, with the finite norm

$$|\mathbf{u}|_{s,p} = \left( \sum_l |u^l|_{s,p}^p \right)^{1/p},$$

we denote by $\mathbb{H}_p^s = \mathbb{H}_p^s(\mathbf{R}^d)$. Similarly, we denote by $\mathbb{H}_p^s(Y) = \mathbb{H}_p^s(\mathbf{R}^d, Y)$ the space of all vector functions $g = (g^l)_{1 \leq l \leq d}$, with $Y$-valued components $g^l$, $1 \leq l \leq d$, so that $\|g\|_{s,p} = (\sum_l |g^l|_{s,p}^p)^{1/p} < \infty$. The set of all infinitely differentiable vector functions $u = (u^1, \ldots, u^d)$ on $\mathbf{R}^d$ with compact support will be denoted by $\mathbb{C}_0^\infty$. We denote by $\mathbb{C}_0^\infty(Y)$ the set of all infinitely differentiable vector functions $u = (u^1, \ldots, u^d)$ on $\mathbf{R}^d$ with compact support (all $u^l$ are $Y$-valued).

When $s = 0$, $\mathbb{H}_p^s(Y) = \mathbb{L}_p(Y) = \mathbb{L}_p(\mathbf{R}^d, Y)$. Also, in this case, the norm $\|g\|_{0,p}$ is denoted more briefly by $\|g\|_p$. To forcefully distinguish $L_p$-norms in spaces of $Y$-valued functions, we write $\|\cdot\|_p$, while in all other cases a norm is denoted by $|\cdot|$.

The duality $\langle \cdot, \cdot \rangle_s$ between $\mathbb{H}_q^s(\mathbf{R}^d)$ and $\mathbb{H}_p^{-s}(\mathbf{R}^d)$, $p \geq 2, s \in (-\infty, \infty)$, and $q = p/(p-1)$ is defined by

$$\langle \phi, \psi \rangle_s = \langle \phi, \psi \rangle_{s,p} = \sum_{i=1}^{d} \int_{\mathbf{R}^d} [\Lambda^s \phi^i](x) \Lambda^{-s} \psi^i(x)\, dx, \qquad \phi \in \mathbb{H}_q^s, \ \psi \in \mathbb{H}_p^{-s}.$$

## 2. Navier–Stokes equation in $\mathbf{R}^d$.

2.1. *Assumptions and main results.* We will consider a stochastic Navier–Stokes equation on $\mathbf{R}^d$ in a finite time interval $[0, T]$. The derivation presented in [27, 29] suggests the following form of this equation for the unknown functions $u = (u^l)_{1 \leq l \leq d}, p, \tilde{p}$:

$$\begin{aligned}
\partial_t u^l(t) &= \partial_i(a^{ij}(t)\,\partial_j u^l(t)) - u^k(t)\,\partial_k u^l(t) - \partial_l p(t) \\
&\quad + b^j(t)\,\partial_j \mathbf{u}(t,x) + \partial_k \tilde{p}(t) h^{l,k}(t) + f^l(t, \mathbf{u}(t)) + \partial_j(f^{l,j}(t, \mathbf{u}(t))) \\
&\quad + [\sigma^k(t)\,\partial_k u^l(t,x) + g^l(t, \mathbf{u}(t)) - \partial_l \tilde{p}(t)] \cdot \dot{W}_t,
\end{aligned}$$
(2.1)
$$\mathbf{u}(0, x) = \mathbf{u}_0(x), \qquad l = 1, \ldots d, \qquad \operatorname{div} \mathbf{u}(t) = 0 \quad \text{in } \mathbf{R}^d \text{ for all } t \in [0, T],$$

where $W$ is a cylindrical Wiener process in a separable Hilbert space $Y$. If $(e_k)$ is a complete orthonormal system (CONS) in $Y$,

$$W_t = \sum_{k=1}^{\infty} W_t^k e_k,$$

where $W_t^k$ are independent standard scalar Wiener processes. In a standard way, for a $Y$-valued adapted random function

$$f_s = \sum_k f_s^k e_k,$$

we define a scalar-valued stochastic integral

$$\int_0^t f_s \cdot dW_s = \sum_k \int_0^t f_s^k\, dW_s^k$$

(in differential form we write $f_s \cdot \dot{W}_s$).

In vector form, we write (2.1) as

$$\begin{aligned}
\partial_t \mathbf{u}(t) &= \partial_i(a^{ij}(t)\,\partial_j \mathbf{u}(t)) - u^k(t)\,\partial_k \mathbf{u}(t) - \nabla p(t) \\
&\quad + b^j(t)\,\partial_j \mathbf{u}(t,x) + \partial_j \tilde{p}(t)\mathbf{h}^j(t) + \mathbf{f}(t,\mathbf{u}(t)) + \partial_j(\mathbf{f}^j(t,\mathbf{u}(t))) \\
&\quad + [\sigma^k(t)\,\partial_k \mathbf{u}(t,x) + \mathbf{g}(t,\mathbf{u}(t)) - \nabla \tilde{p}(t)] \cdot \dot{W}_t,
\end{aligned} \qquad (2.2)$$

$$\mathbf{u}(0) = \mathbf{v}, \qquad \mathrm{div}\,\mathbf{u}(t) = 0 \quad \text{in } \mathbf{R}^d \text{ for all } t \in [0,T],$$

where $\mathbf{h}^j = (h^{l,j})_{1\leq l\leq d}$, $\mathbf{f}^j = (f^{l,j})_{1\leq l\leq d}$, $j = 1,\ldots,d$, $\mathbf{f} = (f^l)_{1\leq l\leq d}$, $\mathbf{g} = (g^l)_{1\leq l\leq d}$. Since $\mathrm{div}\,u(t) = 0$, using the Helmholtz decomposition of vector fields (see the Appendix) we have

$$\begin{aligned}
\nabla \tilde{p}(t,x) &= \tilde{\mathbf{L}}(\mathbf{u}(t),t) = (\tilde{L}_l(\mathbf{u}(t),t))_{1\leq l\leq d} \\
&= \mathcal{G}(\sigma^k(t)\,\partial_k \mathbf{u}(t) + \mathbf{g}(t,\mathbf{u}(t)))
\end{aligned}$$

and

$$\begin{aligned}
\nabla p(t,x) &= \mathcal{G}[-u^k(t)\,\partial_k \mathbf{u}(t) + \partial_i(a^{ij}(t)\,\partial_j \mathbf{u}(t))(t) + \mathbf{f}(t,\mathbf{u}(t)) \\
&\quad + b^i(t)\,\partial_i \mathbf{u}(t) + \tilde{L}_i(\mathbf{u}(t),t)\mathbf{h}^i(t) + \partial_i(\mathbf{f}^i(t,\mathbf{u}(t)))],
\end{aligned}$$

where $\mathcal{G}$ and $\mathcal{S}$ are the projection operators defined in the Appendix.

Thus, instead of (2.1), we can consider the following equivalent equation:

$$\begin{aligned}
\partial_t \mathbf{u}(t) &= \mathcal{S}[\partial_i(a^{ij}(t)\,\partial_j \mathbf{u}(t)) - u^k(t)\,\partial_k \mathbf{u}(t) + b^i(t)\,\partial_i \mathbf{u}(t) \\
&\quad + \tilde{L}_i(\mathbf{u}(t),t)\mathbf{h}^i(t) + \mathbf{f}(t,\mathbf{u}(t)) + \partial_i(\mathbf{f}^i(t,\mathbf{u}(t)))] \\
&\quad + \mathcal{S}[\sigma^i(t)\,\partial_i \mathbf{u}(t) + \mathbf{g}(t,\mathbf{u}(t))]\dot{W}_t,
\end{aligned} \qquad (2.3)$$

$$\mathbf{u}(0) = \mathbf{u}_0, \qquad t \in [0,T].$$

Everywhere below it will be assumed that:

(i) $a^{ij}, b^j$ are measurable functions on $[0,T] \times \mathbf{R}^d$, $f^{l,j}, f^l$ are measurable functions on $[0,T] \times \mathbf{R}^d \times \mathbf{R}^d$;

(ii) $\sigma^j, h^{l,j}$ are $Y$-valued measurable functions on $[0,T] \times \mathbf{R}^d$, $g^l$ are $Y$-valued functions on $[0,T] \times \mathbf{R}^d \times \mathbf{R}^d$; and

(iii) matrix $(a^{ij})$ is nonnegative.

In addition we assume the following:

(B1) $|a^{ij}|, |b^j|, |f^{l,j}|, |\sigma^j|_Y, |\,\mathrm{div}\,\sigma|_Y, |h^{l,j}|_Y$ are bounded by a constant $K$, and there is $\delta > 0$ such that, for all $\xi \in \mathbf{R}^d$,

$$[a^{ij}(t,x) - \tfrac{1}{2}\sigma^j(t,x) \cdot \sigma^i(t,x)]\xi_i \xi_j \geq \delta|\xi|^2;$$

(B2) there exist a constant $C$ and a measurable function $H(t,x)$ on $[0,T] \times \mathbf{R}^d$ such that
$$|f^l(t,x,\mathbf{u})| + |f^{l,j}(t,x,\mathbf{u})| + |g^l(t,x,\mathbf{u})|_Y \leq C|\mathbf{u}| + H(t,x),$$
and for all $t,x$, the functions $f^l(t,x,\mathbf{u}), f^{l,j}(t,x,\mathbf{u}), g^l(t,x,\mathbf{u})$ are continuous in $\mathbf{u}$, where
$$\int_0^T \int_{\mathbf{R}^d} H(t,x)^2 \, dt \, dx < \infty.$$

REMARK 2.1. Note that in (B1) the derivatives $\partial_i \sigma^i$ are understood as Schwartz distributions, but it is assumed that $\operatorname{div}\sigma := \sum_{i=1}^d \partial_i \sigma^i$ is a bounded $Y$-valued function. Obviously, the latter assumption holds in the important case when $\sum_{i=1}^d \partial_i \sigma^i = 0$.

The main results of the paper are the following two statements.

THEOREM 2.1. *Let* (B1) *and* (B2) *hold and let* $u_0 \in \mathbb{L}_2$. *Then there exist a probability space* $(\Omega, \mathcal{F}, \mathbf{P})$ *with a right-continuous filtration* $\mathbb{F} = (\mathcal{F}_t)$ *of $\sigma$-algebras, a cylindrical $\mathbb{F}$-adapted Wiener process $W_t$ in $Y$, and an $\mathbb{L}_2$-valued weakly continuous $\mathbb{F}$-adapted process $u(t)$ such that*
$$\mathbf{E}\left[\sup_{s \leq T} |\mathbf{u}(s)|_2^2 + \int_0^T |\nabla \mathbf{u}(s)|_2^2 \, ds\right] < \infty$$
*and* (2.3) *holds. Moreover, if $d = 2$, then $u(t)$ is (strongly) continuous in $t$.*

THEOREM 2.2. *Let $d = 2$, $u_0 \in \mathbb{L}_2$, let* (B1), (B2) *hold, and for all $l, j, t, x$ and every $\mathbf{u}, \bar{\mathbf{u}}$,*
$$|f^l(t,x,\mathbf{u}) - f^l(t,x,\bar{\mathbf{u}})| + |g^l(t,x,\mathbf{u}) - g^l t,x,\bar{\mathbf{u}})|_Y$$
$$+ |f^{l,j}(t,x,\mathbf{u}) - f^{l,j}(t,x,\bar{\mathbf{u}})|$$
$$\leq K|\mathbf{u} - \bar{\mathbf{u}}|.$$
*Let $(\Omega, \mathcal{F}, \mathbf{P})$ be a probability space with a right-continuous filtration $\mathbb{F} = (\mathcal{F}_t)$ of $\sigma$-algebras and a cylindrical $\mathbb{F}$-adapted Wiener process $W_t$ in $Y$.*

*Then there is a pathwise unique continuous $\mathbb{L}_2$-valued $\mathbb{F}$-adapted solution $u(t)$ to* (2.3) *such that*
$$\mathbf{E}\left[\sup_{s \leq T} |\mathbf{u}(s)|_2^2 + \int_0^T |\nabla \mathbf{u}(s)|_2^2 \, ds\right] < \infty.$$
*Moreover, the distributions of the solutions on different probability spaces coincide.*

In fact in two dimensions, we prove the convergence in probability to $\mathbf{u}(t)$ of the approximating sequence $\mathbf{u}_n(t)$ constructed below.

2.2. *Approximations of Navier–Stokes equations.*

2.2.1. *Construction of approximating sequence.* Let $\psi, \varphi \in C_0^\infty(\mathbf{R}^d)$, $\psi, \varphi \geq 0$, $\int \psi\, dx = \int \varphi\, dx = 1$, $\psi_\varepsilon(x) = \varepsilon^{-d}\psi(x/\varepsilon), \varphi_\varepsilon(x) = \varepsilon^{-d}\varphi(x/\varepsilon)$. Let $\zeta \in C_0^\infty(\mathbf{R}^d)$, $0 \leq \zeta \leq 1$, $\zeta(x) = 1$ if $|x| \leq 1$, $\zeta(x) = 0$ if $|x| > 2$. Using these functions, we mollify the coefficients and the functions of the equation (2.3). Let

$$A^{ij}(t,x) = a^{ij}(t,x) - \tfrac{1}{2}\sigma^j(t,x) \cdot \sigma^i(t,x).$$

Define

$$A_n^{ij}(t,x) = A^{ij}(t,\cdot) * \psi_{1/n}(x) = \int \psi_{1/n}(x-y) A^{ij}(t,y)\, dy,$$

$$\sigma_n^j(t,x) = \sigma^j(t,\cdot) * \psi_{1/n}(x) = \int \psi_{1/n}(x-y) \sigma^j(t,y)\, dy,$$

$$a_n^{ij}(t,x) = A_n^{ij}(t,x) + \tfrac{1}{2}\sigma_n^j(t,x) \cdot \sigma_n^i(t,x).$$

Then,

(2.4) $$[a_n^{ij}(t,x) - \tfrac{1}{2}\sigma_n^j(t,x) \cdot \sigma_n^i(t,x)]\xi_i \xi_j \geq \delta |\xi|^2.$$

Let

$$\tilde{f}_n^l(t,x,\mathbf{u}) = \mathbb{1}_{\{|x|\leq n\}} \mathbb{1}_{\{|H(t,x)|\leq n\}} f^l(t,x,\mathbf{u}) \zeta(\mathbf{u}/n),$$

$$\tilde{f}_n^{l,j}(t,x,\mathbf{u}) = \mathbb{1}_{\{|x|\leq n\}} \mathbb{1}_{\{|H(t,x)|\leq n\}} f^{l,j}(t,x,\mathbf{u}) \zeta(\mathbf{u}/n),$$

$$\tilde{g}_n^l(t,x,\mathbf{u}) = \mathbb{1}_{\{|x|\leq n\}} \mathbb{1}_{\{|H(t,x)|\leq n\}} g^l(t,x,\mathbf{u}) \zeta(\mathbf{u}/n),$$

$$h_n^{l,j}(t,x) = h^{l,j}(t,\cdot) * \psi_{1/n}(x) = \int \psi_{1/n}(x-y) h^{l,j}(t,y)\, dy,$$

$$b_n^j(t,x) = b^{l,j}(t,\cdot) * \psi_{1/n}(x) = \int \psi_{1/n}(x-y) b^{l,j}(t,y)\, dy.$$

Define

$$\tilde{f}_{n,\varepsilon}^l(t,x,\mathbf{u}) = \tilde{f}_n^l(t,x,\cdot) * \varphi_\varepsilon(\mathbf{u}) = \int \tilde{f}_n^l(t,x,\bar{u}) \varphi_\varepsilon(\mathbf{u}-\bar{u})\, d\bar{u},$$

$$\tilde{f}_{n,\varepsilon}^{l,j}(t,x,\mathbf{u}) = \tilde{f}_n^{l,j}(t,x,\cdot) * \varphi_\varepsilon(\mathbf{u}) = \int \tilde{f}_n^{l,j}(t,x,\bar{u}) \varphi_\varepsilon(\mathbf{u}-\bar{u})\, d\bar{u},$$

$$\tilde{g}_{n,\varepsilon}^l(t,x,\mathbf{u}) = \tilde{g}_n^l(t,x,\cdot) * \varphi_\varepsilon(\mathbf{u}) = \int \tilde{g}_n^l(t,x,\bar{u}) \varphi_\varepsilon(\mathbf{u}-\bar{u})\, d\bar{u},$$

and choose $\varepsilon_n \to 0$ so that

$$0 = \lim_n \int_0^T \int \left[ \sup_{\mathbf{u}} (|\tilde{f}_{n,\varepsilon_n}^l(t,x,\mathbf{u}) - \tilde{f}_n^l(t,x,\mathbf{u})|^2 \right.$$

$$\text{(2.5)} \qquad + |\tilde{f}^{l,j}_{n,\varepsilon_n}(t,x,\mathbf{u}) - \tilde{f}^{l,j}_n(t,x,\mathbf{u})|^2$$

$$+ |\tilde{g}^l_{n,\varepsilon_n}(t,x,\mathbf{u}) - \tilde{g}^l_n(t,x,\mathbf{u})|^2) + \varepsilon_n^2 \mathbb{1}_{\{|x| \le n\}} \bigg] dx.$$

Let $f^l_n(t,x,\mathbf{u}) = \tilde{f}^l_{n,\varepsilon_n}(t,x,\mathbf{u})$, $f^{l,j}_n(t,x,\mathbf{u}) = \tilde{f}^{l,j}_{n,\varepsilon_n}(t,x,\mathbf{u})$, $g^l_n(t,x,\mathbf{u}) = \tilde{g}^l_{n,\varepsilon_n}(t,x,\mathbf{u})$ and $\mathbf{f}_n = (f^l_n)_{1 \le l \le d}, \mathbf{g}_n = (g^l_n)_{1 \le l \le d}, \mathbf{f}^j_n = (f^{l,j}_n)_{1 \le l \le d}, \mathbf{h}^j_n = (h^{l,j}_n)_{1 \le l \le d}$.

For each $\mathbf{v} \in \mathbb{H}^1_2$, we define

$$\tilde{\mathbf{L}}_n(\mathbf{v},t) = (\tilde{L}_{n,l}(\mathbf{v},t))_{1 \le l \le d} = \mathcal{G}(\sigma^j_n(t)\,\partial_j\mathbf{v} + \mathbf{g}_n(t,\mathbf{v})),$$

$$A_n(t,\mathbf{v}) = \mathcal{S}[\partial_i(a^{ij}_n(t)\,\partial_j\mathbf{v}) + b^j_n(t)\,\partial_j\mathbf{v} + \tilde{L}_{n,j}(\mathbf{v},t)\mathbf{h}^j(t) + \mathbf{f}_n(t,\mathbf{v}) + \partial_j(\mathbf{f}^j_n(t,\mathbf{v}))],$$

$$B_n(t,\mathbf{v}) = \mathcal{S}[\sigma^i_n(t)\,\partial_i\mathbf{v} + \mathbf{g}_n(t,\mathbf{v})]$$

and

$$N_n(t,\mathbf{v}) = \mathcal{S}[\Psi_n(v^k)\,\partial_k\mathbf{v}],$$

where $\Psi_n(v^k) = \mathbf{v}^k * \psi_{1/n}$.

LEMMA 2.3. *Let* (B1) *and* (B2) *be satisfied. Then there is a constant $C$ independent of $n$ such that for all $\mathbf{v} \in \mathbb{H}^1_2$,*

$$|A_n(t,\mathbf{v})|_{-1,2} \le C[|\mathbf{v}|_{1,2} + |H(t)|_2 + |H_n(t)|_2],$$
$$|B_n(t,\mathbf{v})|_{-1,2} \le C[|\mathbf{v}|_2 + |H(t)|_2 + |H_n(t)|_2],$$

*where $H_n(t) = H_n(t,x)$ is a deterministic function so that $\lim_n \int_0^T |H_n(t)|^2_2\, dt = 0$.*

*Also, for each $k_0 > (d/2) + 1$, there is a constant $C$ independent of $n$ such that for all $\mathbf{v} \in \mathbb{H}^1_2$ so that $\text{div}\,\mathbf{v} = 0$,*

$$|N_n(t,\mathbf{v})|_{-k_0,2} \le C|\mathbf{v}|^2_2.$$

PROOF. For $\mathbf{v} \in \mathbb{H}^1_2$, we have, by Lemma A.1,

$$|\mathcal{S}[\partial_i(a^{ij}_n(t)\,\partial_j\mathbf{v}) + \partial_i(\mathbf{f}^i_n(t,\mathbf{v}))]|_{-1,2}$$
$$\le C[|\mathcal{S}[(a^{ij}_n(t)\,\partial_j\mathbf{v}) + \mathbf{f}^i_n(t,\mathbf{v})]|_2]$$
$$\le [|(a^{ij}_n(t)\,\partial_j v^l) + f^{l,i}_n(t,\mathbf{v})|_2]$$
$$\le C[|\mathbf{v}|_{1,2} + |H(t)|_2 + |H_n(t)|_2].$$

Similarly,
$$|\mathcal{S}[b_n^j(t)\,\partial_j\mathbf{v} + \tilde{L}_{n,i}(\mathbf{v},t)\mathbf{h}^i(t) + \mathbf{f}_n(t,\mathbf{v})]|_{-1,2}$$
$$\leq |\mathcal{S}[b_n^j(t)\,\partial_j\mathbf{v} + \tilde{L}_{n,i}(\mathbf{v},t)\mathbf{h}^i(t) + \mathbf{f}_n(t,\mathbf{v})]|_2$$
$$\leq C|b_n^j(t)\,\partial_j v^l + \tilde{L}_{n,i}(\mathbf{v},t)\mathbf{h}^i(t) + \mathbf{f}_n(t,\mathbf{v})|_2$$
$$\leq C[|\mathbf{v}|_{1,2} + |H(t)|_2 + |H_n(t)|_2].$$

Since $|\partial_i \sigma_n^i(t)| \leq K$, we have by (2.8) and Lemma A.1,
$$|\mathcal{S}[\sigma_n^i(t)\,\partial_i\mathbf{v}]|_{-1,2} \leq |\partial_i(\sigma_n^i(t)\mathbf{v})|_{-1,2} + |\partial_i \sigma_n^i(t)\mathbf{v}|_2 \leq C|\mathbf{v}|_2,$$
$$|\mathcal{S}[\mathbf{g}_n(t,\mathbf{v})]|_2 \leq C[|\mathbf{v}|_2 + |H(t)|_2 + |H_n(t)|_2].$$

Let $k_0 > (d/2) + 1$. Then, obviously,
$$|\mathcal{S}[\Psi_n(v^k)\,\partial_k\mathbf{v}]|_{-k_0,2} \leq |\Psi_n(v^k)\,\partial_k\mathbf{v}|_{-k_0,2}.$$

Since $\Psi_n(v^k)\,\partial_k\mathbf{v} = \partial_k(\Psi_n(v^k)\mathbf{v})$, we have for each $\bar{\mathbf{v}} \in \mathbb{C}_0^\infty \subseteq \mathbb{H}_2^{k_0}$
$$\langle \Psi_n(v^k)\,\partial_k\mathbf{v}, \bar{\mathbf{v}}\rangle_{k_0} = -\int (\Psi_n(v^k)\mathbf{v}, \partial_k\bar{\mathbf{v}})\,dx$$

and by Sobolev's embedding theorem,
$$|\langle \Psi_n(v^k)\,\partial_k\mathbf{v}, \bar{\mathbf{v}}\rangle_{k_0}| \leq \sup_x |\nabla\bar{\mathbf{v}}(x)| \int |\Psi_n(v^k)\mathbf{v}|\,dx$$
$$\leq C|\bar{\mathbf{v}}|_{k_0,2}|\Psi_n(v^k)|_2|\mathbf{v}|_2 \leq C|\bar{\mathbf{v}}|_{k_0,2}|\mathbf{v}|_2^2.$$

So,
$$|\mathcal{S}[\Psi_n(v^k)\,\partial_k\mathbf{v}]|_{-k_0,2} \leq C|\mathbf{v}|_2^2. \qquad \square$$

Now we construct a sequence of approximations. For each $n$, we find $u = (u^l)_{1\leq l\leq d} = u_n = (\mathbf{u}^{n,l})_{1\leq l\leq d}$ by solving

$$\partial_t u^l(t) = \partial_i(a_n^{ij}(t)\,\partial_j u^l(t)) - \Psi_n(u^k(t))\,\partial_k u^l(t)$$
$$- \partial_l p(t) + b_n^j(t)\,\partial_j u^l(t) + \tilde{L}_{n,j}(\mathbf{u}(t),t)h_n^{l,j}(t)$$
(2.6)
$$+ f_n^l(t, \mathbf{u}(t)) + \partial_i(f_n^{l,i}(t,\mathbf{u}(t)))$$
$$+ [\sigma_n^i(t)\,\partial_i u^l(t) + g_n^l(t,\mathbf{u}(t)) - \partial_l \tilde{p}(t)]\dot{W}_t,$$
$$\mathbf{u}(0) = \mathbf{u}_{0,n}, \qquad \mathrm{div}\,\mathbf{u}(t) = 0 \qquad \text{for all } t \in [0,T],$$

where $\Psi_n(u^k(t)) = u^k(t,\cdot) * \psi_{1/n}(x)$, $u_{0,n} = u_0 * \psi_{1/n}$. Equivalently, $u_n(t)$ satisfies

$$\partial_t \mathbf{u}(t) = \mathcal{S}[\partial_i(a_n^{ij}(t)\,\partial_j\mathbf{u}(t)) - \Psi_n(u^k(t))\,\partial_k \mathbf{u}(t)$$
$$+ b_n^j(t)\,\partial_j \mathbf{u}(t) + \tilde{L}_{n,j}(\mathbf{u}(t),t)\mathbf{h}^j(t) + \mathbf{f}_n(t,\mathbf{u}(t)) + \partial_j(\mathbf{f}_n^j(t,\mathbf{u}(t)))]$$
(2.7)
$$+ \mathcal{S}[\sigma_n^j(t)\,\partial_j \mathbf{u}(t) + \mathbf{g}_n(t,\mathbf{u}(t))]\dot{W}_t,$$
$$\mathbf{u}(0) = \mathbf{u}_{0,n} \qquad \text{for all } t.$$

PROPOSITION 2.4. *Let* (B1), (B2) *be satisfied and let* $E|u_0|_2^2 < \infty$. *Then for each* $n$, *there exists a unique* $\mathbb{L}_2$-*valued continuous solution* $u_n(t)$ *of* (2.6) *such that* $\int_0^T |\nabla u_n(s)|_2^2 \, ds < \infty$, $\mathrm{div}\, u_n(t) = 0$ *for all* $t$, **P**-*a.s. Moreover,*

$$\sup_n \mathbf{E}\left[\sup_{t \leq T} |\mathbf{u}_n(t)|_2^2 + \int_0^T |\nabla \mathbf{u}_n(t)|_2^2 \, dt\right] < \infty.$$

PROOF. It is readily checked that for each $n$, there is a constant $K_n$ so that for all $t, x, \mathbf{u}.\bar{\mathbf{u}}, |\alpha| \leq 2$,

$$|\partial^\alpha a_n^{ij}(t,x)| + \|\partial^\alpha \sigma_n^j(t,x)\| + |\partial^\alpha b_n^j(t,x)| \leq K_n,$$
$$|f_n^l(t,x,\mathbf{u}) - f_n^l(t,x,\bar{\mathbf{u}})| + |f_n^{l,j}(t,x,\mathbf{u}) - f_n^{l,j}(t,x,\bar{\mathbf{u}})| \leq K_n|\mathbf{u} - \bar{\mathbf{u}}|,$$
$$|g_n^l(t,x,\mathbf{u}) - g_n^l(t,x,\bar{\mathbf{u}})|_Y \leq K_n|\mathbf{u} - \bar{\mathbf{u}}|$$

and

$$[a_n^{ij}(t,x) - \tfrac{1}{2}(\sigma_n^j(t,x), \sigma_n^i(t,x))_Y]\xi_i\xi_j \geq \tfrac{1}{2}\delta|\xi|^2.$$

Also, there is a constant $C$ independent of $n$ such that

(2.8)
$$|a_n^{ij}(t,x)| + |b_n^j(t,x)| + |\sigma_n^j(t,x)|_Y \leq C,$$
$$|f_n^l(t,x,\mathbf{u})| + |f_n^{l,j}(t,x,\mathbf{u})| + |g_n^l(t,x,\mathbf{u})|_Y \leq H(t,x) + H_n(x) + C|\mathbf{u}|,$$

where $\int_0^T |H(t)|_2^2 \, dt < \infty$, and $H_n(x) = C\varepsilon_n^2 \mathbb{1}_{\{|x| \leq n\}}$. According to (2.5),

(2.9) $$\lim_n |H_n|_2^2 \, dt = 0.$$

For each $\mathbf{v}, \bar{\mathbf{v}} \in \mathbb{H}_2^1, t$,

$$|(\tilde{L}_{n,j}(\mathbf{v},t) - \tilde{L}_{n,j}(\bar{\mathbf{v}},t))\mathbf{h}_n^j(t)|_{-1,2}$$
$$\leq |\partial_i[\mathbf{h}_n^j(t)\mathcal{G}(\sigma_n^i(t)(\mathbf{v} - \bar{\mathbf{v}}))_j]|_{-1,2}$$
$$+ |\partial_i[\mathbf{h}_n^j(t)\mathcal{G}(\sigma_n^i(t)(\mathbf{v} - \bar{\mathbf{v}}))_j]|_2$$
$$+ C|\partial_i\sigma^i(t)(\mathbf{v} - \bar{\mathbf{v}})|_2 + |g_n(t,\mathbf{v}) - g_n(t,\bar{\mathbf{v}})|_2$$
$$\leq K_n|\mathbf{v} - \bar{\mathbf{v}}|_2.$$

Since all the assumptions of Propositions 1 and 3 in [29] are satisfied, there is a unique $\mathbb{L}_2$-valued continuous solution $u_n(t)$ of (2.6) such that $\int_0^T |\nabla u_n(r)|_2^2 \, ds < \infty$, **P**-a.s. Obviously, $\mathrm{div}\, u_n(t) = \mathrm{div}\, \Psi_\kappa(u_n(t)) = 0$ for all $t$, and

$$\int_0^t \int \Psi_n(\mathbf{u}_n)^k(r)\, \partial_k u_n^l(r) u_n^l(r) \, dx \, ds = 0.$$

Therefore, by the Itô formula for $|u_n(t)|_2^2$, we have

$$
\begin{aligned}
|\mathbf{u}_n(t)|_2^2 = |\mathbf{u}_0|_2^2 &+ \int_0^t \int \{2[-a_n^{ij}(s)\, \partial_i u_n^l(s) - f_n^{l,j}(s, \mathbf{u}_n(s))]\, \partial_j u_n^l(s) \\
&+ 2[b_n^{l,k}(s)\, \partial_k u_n^l(s) + f_n^l(s, \mathbf{u}_n(s)) \\
&\qquad + \tilde{L}_n^k(s, \mathbf{u}_n(s)) \cdot h^{l,k}(s)] u_n^l(s) \\
&+ |\mathcal{S}[\sigma_n^j(s)\, \partial_j u_n^l(s) + g_n^l(s, \mathbf{u}_n(s))]|_Y^2\} \, dx\, ds \\
&+ 2 \int_0^t \int [\sigma_n^j(s)\, \partial_j u_n^l(s) + g_n^l(s, \mathbf{u}_n(s))] u_n^l(s) \, dx\, dW_s.
\end{aligned}
\tag{2.10}
$$

Let $\tau$ be an arbitrary stopping time such that

$$
\mathbf{E} \int_0^\tau \left| \int [\sigma^{l,j}(s)\, \partial_j u_n^l(s) + g^l(s, \mathbf{u}_n(s))] u_n^l(s)\, dx \right|_Y^2 ds < \infty.
\tag{2.11}
$$

Since, by Lemma A.1,

$$
\int |\mathcal{S}[\sigma^{l,j}(s)\, \partial_j u_n^l(s) + g^l(s, \mathbf{u}_n(s))]|_Y^2 \, dx
$$

$$
\leq C \int |\sigma^{l,j}(s)\, \partial_j u_n^l(s) + g^l(s, \mathbf{u}_n(s))|_Y^2 \, dx,
$$

using standard arguments (see, e.g., [30], Section 4.1), we obtain that there are constants $\varepsilon, C > 0$ independent of $n$ and $\tau$ such that for all $t$,

$$
\mathbf{E}\left[|\mathbf{u}_n(t \wedge \tau)|_2^2 + \varepsilon \int_0^{\tau \wedge t} |\nabla \mathbf{u}_n(s)|_2^2 \, ds\right]
$$

$$
\leq \mathbf{E}|\mathbf{u}_0|_2^2 + C\mathbf{E} \int_0^{t \wedge \tau} (|H(s)|_2^2 + |\mathbf{u}_n(s)|_2^2) \, ds.
$$

So, by Gronwall's inequality, there is a constant $C$ independent of $n, \tau$ such that

$$
\mathbf{E}\left[|\mathbf{u}_n(t \wedge \tau)|_2^2 + \int_0^{\tau \wedge t} |\partial \mathbf{u}_n(s)|_2^2 \, ds\right] \leq C.
$$

Since $\tau$ is an arbitrary stopping time satisfying (2.11), by Fatou's lemma, for all $t$,

$$
\mathbf{E}\left[|\mathbf{u}_n(t)|_2^2 + \varepsilon \int_0^t |\partial \mathbf{u}_n(s)|_2^2 \, ds\right] \leq C.
\tag{2.12}
$$

Using (2.12), (2.10) and Burkholder's inequality, we easily obtain that

$$
\sup_n \mathbf{E} \sup_{s \leq T} |\mathbf{u}_n(s)|_2^2 < \infty,
$$

and the estimate of the solution follows. □

REMARK 2.2. Note that $u_n(t)$ is a solution of the following equation:

$$\partial_t \mathbf{u}_n(t) = A_n(t, \mathbf{u}_n(t)) - N_n(t, \mathbf{u}_n(t)) + B_n(t, \mathbf{u}_n(t)) \cdot \dot{W}_t,$$

(2.13)
$$\mathbf{u}_n(0) = \mathbf{u}_{0,n}, \qquad \operatorname{div} \mathbf{u}_n(t) = 0.$$

Therefore, combining Proposition 2.4 and Lemma 2.3, we have the following obvious statement.

COROLLARY 2.5. *Let* (B1), (B2) *be satisfied, and let* $E|u_0|_2^2 < \infty$. *Then there is a constant* $C$ *independent of* $n$ *such that* $dt \times d$, $\mathbf{P}$*-a.e.*,

$$|A_n(t, \mathbf{u}_n(t))|_{-1,2} \leq C[|\mathbf{u}_n(t)|_{1,2} + |H(t)|_2 + |H_n|_2],$$

$$|B_n(t, \mathbf{u}(t))|_{-1,2} \leq C[|\mathbf{u}_n(t)|_2 + |H(t)|_2 + |H_n|_2],$$

*where* $H_n = H_n(x)$ *is a deterministic function so that* $\lim_n |H_n|_2 = 0$.

*Also, for each* $k_0 > (d/2) + 1$, *there is a constant* $C$ *independent of* $n$ *such that* $dt \times d$, $\mathbf{P}$*-a.e.*,

$$|N_n(t, \mathbf{u}_n(t))|_{-k_0,2} \leq C|\mathbf{u}_n(t)|_2^2.$$

2.2.2. *Weak compactness of approximations.* For each $n$, the solution $u_n$ of (2.6) induces a measure $\mathbf{P}^n$ on some trajectory space determined by the estimates of Proposition 2.4.

Denote by $\mathbb{L}_{2,\mathrm{loc}}$ the space $\mathbb{L}_2$ with a topology of $L_2$-convergence on compact subsets of $\mathbf{R}^d$. It is defined by the seminorms

$$|\mathbf{v}|_{2;R} = \int_{|x| \leq R} |\mathbf{v}|^2 \, dx, \qquad R > 0.$$

Fix $U = \mathbb{H}_2^{k_0}, k_0 > (d/2) + 1$. Denote by $U'_{\mathrm{loc}}$ the space $U'$ with a topology defined by the seminorms

$$|\mathbf{g}|_{U',R} = \sup\{|\mathbf{g}(\mathbf{v})| : \mathbf{v} \in \mathbb{C}_0^\infty, |\mathbf{v}|_U \leq 1, \operatorname{supp} \mathbf{v} \subset B_R\}, \qquad 0 < R < \infty,$$

where $B_R = \{x : |x| < R\}$.

LEMMA 2.6. *The embedding* $\mathbb{L}_2 \to U'_{\mathrm{loc}}$ *is compact.*

PROOF. Let $\{\mathbf{x}_k\}_{k \geq 1}$ be a bounded set in $\mathbb{L}_2$. Then there exist a subset $\{\mathbf{x}_{k'}\}_{k' \geq 1} \subset \{\mathbf{x}_k\}_{k \geq 1}$ and $\mathbf{x} \in \mathbb{L}_2$ such that

(2.14) $$\lim_{k' \to \infty} (\mathbf{x}_{k'}, f)_2 = (\mathbf{x}, f)_2 \qquad \text{for any } f \in \mathbb{L}_2.$$

Let $\{\mathbf{e}_i\}_{i\geq 1}$ be an orthonormal basis in $\mathbb{H}_2^{-k_0}(B_R)$. Obviously,

$$\sum_{i\geq N}(\mathbf{x}_{k'} - \mathbf{x}, \mathbf{e}_i)^2_{U',R} \leq |\mathbf{x}_{k'} - \mathbf{x}|^2_{\mathbb{H}^{-k_0}(B_R)} \sum_{i\geq N} |\mathbf{e}_i|^2_{\mathbb{H}^{-k_0}(B_R)}$$

(2.15)

$$\leq C \sum_{i\geq N} |(-\Delta_0)^{-k_0/2} \mathbf{e}_i|^2_{\mathbb{H}^{-k_0}(B_R)},$$

where $\Delta_0$ is the Laplace operator on $L_2(B_R)$ with zero boundary conditions. Since $(-\Delta_0)^{-k_0/2}$ is a Hilbert–Schmidt operator, then for any $\varepsilon > 0$, there is $N_\varepsilon$ such that the right-hand side of (2.15) is less than $\varepsilon$ for all $N \geq N_\varepsilon$. Now, compactness follows by (2.14). $\square$

We remark that the lemma holds also for arbitrary $k_0 > 2$. This could be proved using arguments similar to those in Remark III.3.2 in [32].

Let $C_{[0,T]}(U'_{\text{loc}})$ be the set of $U'_{\text{loc}}$-valued trajectories with the topology $\mathcal{T}_1$ of the uniform convergence on $[0,T]$. Let $C_{[0,T]}(\mathbb{L}_{2,w})$ be the set of $\mathbb{L}_2$-valued weakly continuous functions with the topology $\mathcal{T}_0$ of the uniform weak convergence on $[0,T]$. Let $\mathbb{L}_{2,w}(0,T;\mathbb{H}_2^1)$ be the set of $\mathbb{H}_2^1$-valued square integrable functions $\mathbf{f}_s$ on $[0,T]$ with a topology $\mathcal{T}_2$ of weak convergence on finite intervals, that is, the topology defined by the maps

$$\mathbf{f}_s \to \int_0^T \langle \mathbf{f}_s, \mathbf{g}_s \rangle_1 \, ds,$$

where $\mathbf{g}$ is $\mathbb{H}_2^{-1}$-valued such that $\int_0^T |g_s|^2_{-1,2} \, ds < \infty$. Let $\mathbb{L}_2(0,T;\mathbb{L}_{2,\text{loc}})$ be the space of square integrable functions with a topology $\mathcal{T}_3$ generated by seminorms

$$|\mathbf{u}|_{2;T,R} = \int_0^T \int_{|x|\leq R} |\mathbf{u}(t,x)|^2 \, dx \, dt, \qquad R > 0.$$

LEMMA 2.7 (cf. [25, 33]). Let $\mathcal{Z} = C_{[0,T]}(U'_{\text{loc}}) \cap C_{[0,T]}(\mathbb{L}_{2,w}) \cap \mathbb{L}_{2,w}(0,T;\mathbb{H}_2^1) \cap \mathbb{L}_2(0,T;\mathbb{L}_{2,\text{loc}})$ and let $\mathcal{T}$ be the supremum of the corresponding topologies. Then $K \subset \mathcal{Z}$ is $\mathcal{T}$-relatively compact if the following conditions hold:

(a) $\sup_{x \in K} \sup_{s \leq T} |x_s|_2 < \infty$,
(b) $\sup_{x \in K} \int_0^T |x_s|^2_{1,2} \, ds < \infty$,
(c) $\lim_{\delta \to 0} \sup_{x \in K} \sup_{|t-s|\leq \delta, \ s,t \leq T} |x_t - x_s|_{U'} = 0$.

PROOF. It can be assumed that $K$ is closed in $\mathcal{T}$. The topologies $\mathcal{T}_0$, $\mathcal{T}_1$, $\mathcal{T}_2$, $\mathcal{T}_3$ are metrizable on $K$. Consider a sequence $(\mathbf{x}^n)$ in $K$. Obviously, (c) yields that $K$ is compact in $\mathcal{T}_2$ topology. By Lemma 2.6, the imbedding $\mathbb{L}_2 \subset U'_{\text{loc}}$ is compact. Therefore, by (a), (c) and the Arzelá–Ascoli theorem for functions taking values in a Fréchet space, there exist a subsequence $(\mathbf{x}^{n_k})$ and a function $\mathbf{x}$ such that $\mathbf{x}^{n_k} \to \mathbf{x}$ in $C_{[0,T]}(U'_{\text{loc}}) \cap \mathbb{L}_{2,w}(0, T; \mathbb{H}_2^1)$ with respect to the supremum of $\mathcal{T}_1$ and $\mathcal{T}_2$.

Since, for all $u \in \mathbb{C}_0^\infty$, $(\mathbf{x}_s^{n_k}, \mathbf{u})_{\mathbb{L}_2} = \langle \mathbf{x}_s^{n_k}, \mathbf{u}\rangle \longrightarrow \langle \mathbf{x}_s, \mathbf{u}\rangle$, we have

$$\sup_{s \leq T} |\mathbf{x}_s|_2 = \sup_{s \leq T} \sup_{\mathbf{u} \in \mathbb{C}_0^\infty} \langle \mathbf{x}_s, \mathbf{u}\rangle \Big/ |\mathbf{u}|_2 \leq \sup_{s \leq T} \liminf_{k \to \infty} |\mathbf{x}_s^{n_k}|_2 < \infty.$$

Next, for every $u \in \mathbb{C}_0^\infty$ there exists $R < \infty$ so that

$$\sup_{s \leq T} |(\mathbf{x}_s - \mathbf{x}_s^{n_k}, \mathbf{u})_{\mathbb{L}_2}| \leq \sup_{s \leq T} |\mathbf{x}_s - \mathbf{x}_s^{n_k}|_{U', R} |\mathbf{u}|_U.$$

It is readily checked now that for every $\mathbf{v} \in U$,

$$\lim_{k \to \infty} \sup_{s \leq T} |(\mathbf{x}_s - \mathbf{x}_s^{n_k}, \mathbf{v})_2| = 0.$$

Thus $\mathbf{x}^{n_k} \to \mathbf{x}$ in $C_{[0,T]}(\mathbb{L}_{2,w})$ as well. Obviously, for each $p > 0$, $R > 0$,

$$(2.16) \qquad \int_0^T |\mathbf{x}_s^{n_k} - \mathbf{x}_s|_{U', R}^p \, ds \to 0$$

as $k \to \infty$. We claim that for each $\epsilon > 0$, $R$, there is a constant $C = C_{\epsilon, R}$ such that for all $u \in \mathbb{H}_2^1$,

$$(2.17) \qquad |\mathbf{u}|_{2;R}^2 \leq \epsilon |\mathbf{u}|_{1,2}^2 + C|\mathbf{u}|_{U', R}^2.$$

Indeed, if (2.17) does not hold, there exists $\epsilon > 0$, $R > 0$ and a sequence $u_n \in U$ such that

$$|\mathbf{u}_n|_{2;R}^2 > \epsilon |\mathbf{u}_n|_{1,2}^2 + n|\mathbf{u}_n|_{U';R}^2.$$

Then for $\mathbf{v}_n = |\mathbf{u}_n|_{2;R}^{-1} u_n$ we have

$$1 > \epsilon |\mathbf{v}_n|_{1,2}^2 + n|\mathbf{v}_n|_{U';R}^2.$$

Thus $(\mathbf{v}_n)$ is a bounded sequence in $\mathbb{H}_2^1$ and $|\mathbf{v}_n|_{U';R} \to 0$ as $n \to \infty$. Since the embedding $\mathbb{H}_2^1 \to \mathbb{L}_{2,\text{loc}}$ is compact, $|\mathbf{v}_n|_{2;R} \to 0$. On the other hand, $|\mathbf{v}_n|_{2;R} = 1$ for all $n$, and we have a contradiction. Thus (2.17) holds. Since

$$\sup_n \int_0^T (|\mathbf{x}_s^n|_{1,2}^2 + |\mathbf{x}_s|_{1,2}^2) \, ds < \infty,$$

it follows by (2.16), (2.17) that

$$|\mathbf{x}^{n_k} - \mathbf{x}|_{2;T,R} \to 0. \qquad \square$$

Let $X(t) = X(t, x) = X(t, w) = w(t) = w(t, x)$, $w \in C_{[0,T]}(U'_{\text{loc}})$. Let $\mathcal{D}_t = \sigma(X(s), s \leq t)$, $\mathbb{D} = (\mathcal{D}_{t+})_{0 \leq t < T}$, $\mathcal{D} = \mathcal{D}_T$.

For each $n$, the solution $u_n$ to (2.6) defines a measure $\mathbf{P}^n$ on $(C_{[0,T]}(U'_{\text{loc}}), \mathcal{D})$.

COROLLARY 2.8. *The set $\{\mathbf{P}^n, n \geq 1\}$ is relatively weakly compact on $(\mathcal{Z}, \mathcal{T})$.*

PROOF. By Remark 2.2, $\mathbf{P}^n$-a.s.,

$$dX_t = [A_n(t, X(t)) - N_n(t, X(t))]\,dt + dM^n_t,$$
$$X(0, x) = \mathbf{u}_0(x),$$

where $M^n_t$ is $\mathbb{H}_2^{-1} \subseteq U'$-valued martingale such that for each $\mathbf{v} \in U$,

$$\langle M^n_t, \mathbf{v}\rangle_1^2 - \int_0^t \langle B_n(s, X(s)), \mathbf{v}\rangle_1^2\,ds \in \mathcal{M}_{\mathrm{loc}}(\mathbb{D}, \mathbf{P}^n).$$

Let $\tau_n = \tau_n(X)$ be a sequence of stopping times such that $\tau_n \leq T$. Let $\delta_n$ be a sequence of numbers so that $1 > \delta_n \downarrow 0$. We have by Corollary 2.5,

$$\mathbf{P}^n |M^n_{\tau_n+\delta_n} - M^n_{\tau_n}|^2_{-1,2} = \mathbf{E}\left|\int_{\tau_n(\mathbf{u}_n)}^{\tau_n(\mathbf{u}_n)+\delta_n} B_n(s, \mathbf{u}_n(s)) \cdot dW_s\right|^2_{-1,2}$$

$$\leq C\mathbf{E}\left[\delta_n \sup_{s \leq T} |\mathbf{u}_n(s)|_2^2 + \int_{\tau_n(\mathbf{u}_n)}^{\tau_n(\mathbf{u}_n)+\delta_n} |H(t)|_2^2\,dt\right]$$
$$+ T|H_n|_2^2.$$

Here and below, with a slight abuse of notation, we write $\mathbf{P}^n f$ for an integral of a measurable function $f$ with respect to the measure $\mathbf{P}^n$.

So,

(2.18) $$\lim_n \mathbf{P}^n |M^n_{\tau_n+\delta_n} - M^n_{\tau_n}|^2_{-1,2} = 0.$$

By Corollary 2.5 and the Hölder inequality,

$$\mathbf{P}^n \int_{\tau_n}^{\tau_n+\delta_n} |A_n(t, X(t))|_{-1,2}\,dt$$
$$= \mathbf{E}\int_{\tau_n(\mathbf{u}_n)}^{\tau_n(\mathbf{u}_n)+\delta_n} |A_n(t, \mathbf{u}_n(t))|_{-1,2}\,dt$$
$$\leq C\mathbf{E}\left[\delta_n^{1/2}\left(\int_0^T |\mathbf{u}_n(t)|_{1,2}^2\,dt\right)^{1/2}\right.$$
$$\left.+ \delta_n^{1/2}\left(\int_0^T (|H(t)|_2 + |H_n|_2)^2\,dt\right)^{1/2} dt\right].$$

Therefore,

(2.19) $$\lim_n \mathbf{P}^n \int_{\tau_n}^{\tau_n+\delta_n} |A_n(t, X(t))|_{-1,2}\,dt = 0.$$

Also, by Corollary 2.5,

$$\mathbf{P}^n \int_{\tau_n}^{\tau_n+\delta_n} |N_n(t, X(t))|_{U'}\, dt = \mathbf{E} \int_{\tau_n(\mathbf{u}_n)}^{\tau_n(\mathbf{u}_n)+\delta_n} |N_n(t, \mathbf{u}_n(t))|_{U'}\, dt$$

$$\leq C\delta_n \mathbf{E} \sup_{s \leq T} |\mathbf{u}_n(s)|_2^2.$$

This and (2.18), (2.19) imply that

(2.20) $$\lim_n \mathbf{P}^n |X_{\tau_n+\delta_n} - X_{\tau_n}|_{U'} = 0.$$

Let $\mathbf{P}_t^n$ be the natural restriction of $\mathbf{P}_n$ to $\sigma(X(t))$. By Lemma 2.6, Proposition 2.4 and Prokhorov's theorem for Fréchet spaces (see [1]), the family of measures $\{\mathbf{P}_t^n, n \geq 1\}$ is relatively compact on $U'_{\text{loc}}$. Also, by the Aldous criterion, (2.8) yields that for each $T > 0$, $\eta > 0$,

$$\lim_{\delta \to 0} \sup_n \mathbf{P}^n \left( \sup_{\substack{|s-t| \leq \delta \\ s, t \leq T}} |X_t - X_s|_{U'} > \eta \right) = 0.$$

Therefore, the relative compactness of measures $\{\mathbf{P}^n, n \geq 1\}$ on $\mathcal{Z}$ with supremum topology $\mathcal{T}$ follows by Lemma 2.7 in a standard way (cf. [25, 33]). □

2.2.3. $\mathbf{P}^n$ *as a solution of a martingale problem.* For $\mathbf{v} \in \mathbb{C}_0^\infty(\mathbf{R}^d)$, denote

$$\varphi_n^{\mathbf{v}}(s, X) = i_0 \langle A_n(s, X_s) - N_n(s, X_s), \mathbf{v} \rangle_{k_0} - \tfrac{1}{2} |\langle B_n(s, X_s), \mathbf{v} \rangle_{k_0, Y}|_Y^2,$$

where $i_0^2 = -1$. (Let us recall $U = \mathbb{H}_2^{k_0}, k_0 > d/2$.) Notice that

$$\langle A_n(s, X_s), \mathbf{v} \rangle_{k_0} = -\int (a_n^{ij}(s)\, \partial_j X_s + \mathbf{f}_n^i(s, X_s), \mathcal{S}(\partial_i \mathbf{v}))\, dx$$

$$+ \int (b_n^j(s)\, \partial_j X_s + \tilde{L}_{n,j}(X_s, s) \cdot \mathbf{h}^j(s) + \mathbf{f}_n(s, X_s), \mathcal{S}(\mathbf{v}))\, dx,$$

$$\langle N_n(s, X_s), \mathbf{v} \rangle_{k_0} = \langle \Psi_n(X_s^k)\, \partial_k X_s, \mathcal{S}(\mathbf{v}) \rangle_{k_0} = -\int (\Psi_n(X_s^k) X_s, \mathcal{S}(\partial_k \mathbf{v}))\, dx,$$

$$\langle B_n(s, X_s), \mathbf{v} \rangle_{k_0, Y} = \int (\sigma_n^i(s)\, \partial_i X_s + \mathbf{g}_n(s, X_s), \mathcal{S}(\mathbf{v}))\, dx.$$

Applying the Itô formula to the scalar semimartingale $\langle X_t, \mathbf{v} \rangle = \langle X_t, \mathbf{v} \rangle_0 = \int X_t^l v^l\, dx$, we obtain the following obvious statement.

LEMMA 2.9. *For each $n$, $\mathbf{P}^n$ is a measure on $\mathcal{Z}$ such that for each test function $\mathbf{v} \in \mathbb{C}_0^\infty(\mathbf{R}^d)$,*

$$L_t^{n, \mathbf{v}} = e^{i_0 \langle X_t, \mathbf{v} \rangle} - \int_0^t e^{i_0 \langle X(s), \mathbf{v} \rangle} \varphi_n^{\mathbf{v}}(s, X)\, ds \in \mathcal{M}_{\text{loc}}^c(\mathbb{D}, \mathbf{P}^n).$$

[*We say $\mathbf{P}^n$ is a solution of the martingale problem $(u_0, A_n, B_n)$.*]

For $\mathbf{v} \in \mathbb{H}_2^1$, we set
$$A(t,\mathbf{v}) = \mathcal{S}[\partial_i(a^{ij}(t)\,\partial_j \mathbf{v}) + b^j(t)\,\partial_j \mathbf{v} + \tilde{L}_j(\mathbf{v},t)\mathbf{h}^j(t)$$
$$+ \mathbf{f}(t,\mathbf{v}) + \partial_j(\mathbf{f}^j(t,\mathbf{v}))],$$
$$N(\mathbf{v}) = \mathcal{S}[v^k\,\partial_k \mathbf{v}]$$

and
$$B(t,\mathbf{v}) = \mathcal{S}[\sigma^i(t)\,\partial_i \mathbf{v} + \mathbf{g}(t,\mathbf{v})],$$

where
$$\tilde{\mathbf{L}}(\mathbf{v},t) = (\tilde{L}_j(\mathbf{v},t))_{1 \leq j \leq d} = \mathcal{G}[\sigma^i(t)\,\partial_i \mathbf{v} + \mathbf{g}(t,\mathbf{v})].$$

Let
$$\varphi^{\mathbf{v}}(s,X_s) = i_0 \langle A(s,X_s) - N(s,X_s), \mathbf{v}\rangle_{k_0} - \tfrac{1}{2}|\langle B(s,X_s), \mathbf{v}\rangle_{k_0,Y}|_Y^2.$$

DEFINITION 2.1. We say a probability measure $\mathbf{P}$ on $\mathcal{Z}$ is a solution of the martingale problem $(u_0, A, B)$ if for each $\mathbf{v} \in \mathbb{C}_0^\infty(\mathbf{R}^d)$,
$$L_t^{\mathbf{v}} = e^{\iota_0 \langle X_t, \mathbf{v}\rangle} - \int_0^t e^{\iota_0 \langle X(s), \mathbf{v}\rangle} \varphi^{\mathbf{v}}(s, X_s)\,ds \in \mathcal{M}_{\mathrm{loc}}^c(\mathbb{D}, \mathbf{P}),$$
and $X_0 = u_0$, $\mathbf{P}$-a.s.

2.3. *Existence of weak global solutions.* In a standard way, we obtain the following statement.

THEOREM 2.10. *Assume* (B1) *and* (B2) *are satisfied. Then for each* $u_0 \in \mathbb{L}_2$ *there is a measure* $\mathbf{P}$ *on* $\mathcal{Z}$ *solving the martingale problem* $(u_0, A, B)$ *such that*

(2.21) $$\mathbf{P}\left[\sup_{t \leq T}|X(t)|_2^2 + \int_0^T |X(s)|_{1,2}^2\,ds\right] < \infty.$$

*Moreover, if* $d = 2$, *then* $\mathbf{P}$-*a.s.*
$$\int_0^T |A(s,X_s) - N(s,X_s)|_{-1,2}^2\,ds < \infty.$$

PROOF. We follow the lines of the proof in [25]. Since the set $\{\mathbf{P}^n, n \geq 1\}$ is relatively compact, we can assume that a sequence of measures $(\mathbf{P}^n)$ converges weakly to some measure $\mathbf{P}$ on $\mathcal{Z}$. Let $\omega_n \to \omega$ in $\mathcal{Z}$. Then, by Lemma 2.7,

(2.22) $$\sup_n \left\{\sup_{s \leq T}(|\omega_n(s)|_2 + |\omega(s)|_2) + \int_0^T (|\omega_n(s)|_{1,2}^2 + |\omega(s)|_{1,2}^2)\,ds\right\} < \infty,$$

and for each $R > 0, \mathbf{v} \in \mathbb{C}_0^\infty$

$$\begin{aligned}(2.23)\quad &\int_0^T \int_{|x|\leq R} |\omega_n(s,x) - \omega(s,x)|^2 \, dx \, ds \\ &+ \sup_{s\leq T} \left| \int (\omega_n(s,x) - \omega(s,x), \mathbf{v}(x)) \, dx \right| \to 0,\end{aligned}$$

as $n \to \infty$. Also, for each $\mathbf{w} \in \mathbb{L}_2([0,T]; \mathbb{H}_2^{-1})$,

$$(2.24) \quad \int_0^T \langle \mathbf{w}(s), \omega_n(s) - \omega(s) \rangle_1 \, ds \to 0,$$

as $n \to \infty$.

It follows from (2.22) and (2.23) that the sequence $\omega_n(t)$ is weakly relatively compact in $\mathbb{L}_2([0,T]; \mathbb{H}_2^1)$. This and assumptions (B1) and (B2) imply that the sequence $(L_t^{n,\mathbf{v}}(\omega_n))$ is equicontinuous in $t$ with respect to $n$. Indeed, by Lemma 2.3, there exists a constant $C$ independent of $n$ so that $dt$-a.e.,

$$|N_n(\omega_n(t),t)|_{-k_0,2} \leq C|\omega_m(t)|_2^2,$$
$$|A_n(t,\omega_n(t))|_{-1,2} \leq C[|\omega_n(t)|_{1,2} + |H(t)|_2 + |H_n|_2],$$
$$\|B_n(t,\omega_n(t))\|_{-1,2} \leq C[|\omega_n(t)|_2 + |H(t)|_2 + |H_n|_2].$$

So, $dt$-a.e.

$$\begin{aligned}(2.25)\quad |\varphi_n^{\mathbf{v}}(t,\omega_n(t))| &\leq C(|\mathbf{v}|_{k_0,2} + |\mathbf{v}|_{1,2}^2) \\ &\times \left( |\omega_n(t)|_{1,2} + \sup_{s\leq T} |\omega_n(s)|_2^2 + |H(t)|_2 + |H_n|_2 \right).\end{aligned}$$

Therefore, by the Hölder inequality, there is a constant independent of $n$ such that for each $r < s$,

$$\left| \int_r^s e^{\iota_0 \langle \omega_n(t), \mathbf{v}\rangle} \varphi^{\mathbf{v}}(t,\omega_n(t)) \, dt \right| \leq C(|r-s| + |r-s|^{1/2}),$$

and equicontinuity in $t$ of the sequence $(L_t^{n,\mathbf{v}}(\omega_n))$ follows.

Now we prove that for each $t \in [0,T]$,

$$(2.26) \quad L_t^{n,\mathbf{v}}(\omega_n) \to L_t^{\mathbf{v}}(\omega).$$

By (2.22), (2.23) and (2.25),

$$\sup_{s\leq T} |e^{\iota_0 \langle \omega_n(s),\mathbf{v}\rangle} - e^{\iota_0 \langle \omega(s),\mathbf{v}\rangle}| \to 0,$$

$$\int_0^T |e^{\iota_0 \langle \omega_n(t),\mathbf{v}\rangle} \varphi_n^{\mathbf{v}}(t,\omega_n(t)) \, dt - e^{\iota_0 \langle \omega(t),\mathbf{v}\rangle} \varphi_n^{\mathbf{v}}(t,\omega_n(t))| \, dt \to 0,$$

as $n \to \infty$. Also, notice

$$
\int \tilde{L}_{n,j}(\omega_n(t), t) \cdot (\mathbf{h}_n^j(t), \mathcal{S}(\mathbf{v})) \, dx
$$
(2.27)
$$
= \int [\sigma_n^j(t) \, \partial_j \omega_n(t) + \mathbf{g}_n(t, \omega_n(t))] \cdot \mathcal{G}[(h_n^{j,l}(t) \mathcal{S}(\mathbf{v})^l)_j] \, dx.
$$

Since

$$
0 = \lim_n \int_0^T (|(a_n^{ij}(t) - a^{ij}(t))\mathcal{S}(\partial_j \mathbf{v})|_2^2 + |(b_n^j(t) - b^j(t))\mathcal{S}(\mathbf{v})|_2^2
$$
$$
+ |\sigma_n^j(t) \cdot \mathcal{G}[(h_n^{j,l}(t)\mathcal{S}(\mathbf{v})^l)_j] - \sigma^j(t) \cdot \mathcal{G}[(h^{j,l}(t)\mathcal{S}(\mathbf{v})^l)_j]|_2^2) \, dt,
$$

it follows by (2.22) and (2.24) that

$$
\lim_n \int_0^t e^{\iota_0 \langle \omega(s), \mathbf{v} \rangle} \int [-a_n^{ij}(s) \, \partial_j \omega_n^l(s) \mathcal{S}(\partial_j v)^l + b_n^j(s) \, \partial_j \omega_n^l(s) \mathcal{S}(v)^l] \, dx \, ds
$$
$$
= \int_0^t e^{\iota_0 \langle \omega(s), \mathbf{v} \rangle} \int [-a^{ij}(s) \, \partial_j \omega^l(s) \mathcal{S}(\partial_j v)^l + b^j(s) \, \partial_j \omega^l(s) \mathcal{S}(v)^l] \, dx
$$

and

$$
\lim_n \int_0^t e^{\iota_0 \langle \omega(s), \mathbf{v} \rangle} \int \sigma_n^j(s) \, \partial_j \omega_n(s) \cdot \mathcal{G}[(h_n^{j,l}(s)\mathcal{S}(\mathbf{v})^l)_j] \, dx \, ds
$$
$$
= \int_0^t e^{\iota_0 \langle \omega(s), \mathbf{v} \rangle} \int \sigma^j(s) \, \partial_j \omega(s) \cdot \mathcal{G}[(h^{j,l}(s)\mathcal{S}(\mathbf{v})^l)_j] \, dx \, ds.
$$

By (2.5) and (2.23), it follows for each $m > 0$,

$$
\limsup_n \int_0^T \int |(\mathbf{f}_n(t, \omega_n(t)) - \mathbf{f}(t, \omega(t)), \mathcal{S}(\mathbf{v}))| \, dx \, dt
$$
$$
\leq C \left[ \limsup_n \int_0^T \left( \int_{|x| \leq m} |\mathbf{f}_n(t, \omega_n(t)) - \mathbf{f}(t, \omega(t))|^2 \, dx \right)^{1/2} \right.
$$
$$
\left. + \limsup_n \int_0^T (|\mathbf{f}_n(t, \omega_n(t))|_2 + |\mathbf{f}(t, \omega(t))|_2) \left( \int_{|x| > m} |\mathcal{S}(\mathbf{v})|^2 \, dx \right)^{1/2} \right]
$$
$$
\leq C \left( \int_{|x| > m} |\mathcal{S}(\mathbf{v})|^2 \, dx \right)^{1/2}.
$$

Since $m$ is arbitrarily large, for each $t$,

$$
\lim_n \int_0^t \int e^{\iota_0 \langle \omega(s), \mathbf{v} \rangle} (\mathbf{f}_n(t, \omega_n(s)), \mathcal{S}(\mathbf{v})) \, dx
$$
(2.28)
$$
= \int_0^t \int e^{\iota_0 \langle \omega(s), \mathbf{v} \rangle} (\mathbf{f}(t, \omega(s)), \mathcal{S}(\mathbf{v})) \, dx.
$$

Similarly,

$$\lim_n \int_0^t \int e^{\iota_0 \langle \omega(s), \mathbf{v} \rangle} (\mathbf{f}_n^i(t, \omega_n(s)), \partial_i \mathcal{S}(\mathbf{v}))\, dx$$
$$= \int_0^t \int e^{\iota_0 \langle \omega(s), \mathbf{v} \rangle} (\mathbf{f}^i(t, \omega(s)), \partial_i \mathcal{S}(\mathbf{v}))\, dx$$

and

$$\lim_n \int_0^t e^{\iota_0 \langle \omega(s), \mathbf{v} \rangle} \int \mathbf{g}_n(s, \omega_n(s)) \cdot \mathcal{G}[(h_n^{j,l}(s)\mathcal{S}(\mathbf{v})^l)_j]\, dx\, ds$$
$$= \int_0^t e^{\iota_0 \langle \omega(s), \mathbf{v} \rangle} \int \mathbf{g}(s, \omega(s)) \cdot \mathcal{G}[(h^{j,l}(s)\mathcal{S}(\mathbf{v})^l)_j]\, dx\, ds.$$

Therefore, for each $t$,

$$\lim_n \int_0^t e^{\iota_0 \langle \omega_n(s), \mathbf{v} \rangle} \langle A_n(s, \omega_n(s)), \mathbf{v} \rangle_{k_0}\, ds = \int_0^t e^{\iota_0 \langle \omega(s), \mathbf{v} \rangle} \langle A(s, \omega(s)), \mathbf{v} \rangle_{k_0}\, ds.$$

Since

$$Q_n(s) = \int (\sigma_n^j(s)\, \partial_j \omega(s) + \mathbf{g}_n(s, \omega(s))) \cdot \mathcal{S}(\mathbf{v})\, dx$$
$$= \int [(\mathbf{g}_n(s, \omega_n(s)) - \partial_j \sigma_n^j(s)\omega_n(s)) \cdot \mathcal{S}(\mathbf{v})$$
$$- \sigma_n^j(s)\omega_n(s) \cdot \partial_j \mathcal{S}(\mathbf{v})]\, dx,$$

by (2.22) and (2.23), we see that there is a constant $C$ independent of $n$ such that

$$|Q_n(s)| \leq C(1 + |H(s)|_2).$$

As in the case of (2.28), we obtain

$$\lim_n \int_0^t e^{i_0 \langle \omega(s), \mathbf{v} \rangle} |Q_n(s)|_Y^2\, ds$$
$$= \lim_n \int_0^t e^{i_0 \langle \omega(s), \mathbf{v} \rangle} |\langle B_n(s, \omega_n(s)), \mathbf{v} \rangle_{k_0, Y}|_Y^2\, ds$$
$$= \int_0^t e^{i_0 \langle \omega(s), \mathbf{v} \rangle} \left| \int (\sigma^j(s)\, \partial_j \omega(s) + \mathbf{g}(s, \omega(s))) \cdot \mathcal{S}(\mathbf{v})\, dx \right|_Y^2 ds$$
$$= \int_0^t e^{i_0 \langle \omega(s), \mathbf{v} \rangle} |\langle B(s, \omega(s)), \mathbf{v} \rangle_{k_0, Y}|_Y^2\, ds.$$

Finally, we prove that for each $t$,

(2.29)
$$\lim_n \int_0^t e^{\iota_0 \langle \omega(s), \mathbf{v} \rangle} \langle N_n(s, \omega_n(s)), \mathbf{v} \rangle_{k_0}\, ds$$
$$= \int_0^t e^{\iota_0 \langle \omega(s), \mathbf{v} \rangle} \langle N(s, \omega(s)), \mathbf{v} \rangle_{k_0}\, ds.$$

Let $\zeta \in C_0^\infty(\mathbf{R}^d)$, $0 \leq \zeta \leq 1$, $\zeta(x) = 1$ if $|x| \leq 1$, $\zeta(x) = 0$ if $|x| > 2$. First of all, for each $m > 0$,

$$\lim_n \int_0^t e^{\iota_0 \langle \omega(s), \mathbf{v} \rangle} \int (\Psi_n(\omega_n^k(s))\omega_n(s), \mathcal{S}(\partial_k \mathbf{v}))\zeta(x/m)\, dx \tag{2.30}$$
$$= \int_0^t e^{\iota_0 \langle \omega(s), \mathbf{v} \rangle} \int (\omega^k(s)\omega(s), \mathcal{S}(\partial_k \mathbf{v}))\zeta(x/m)\, dx.$$

Indeed,

$$\int_0^T \int |((\Psi_n(\omega_n^k(s)) - \omega_n^k(s))\omega_n(s), \mathcal{S}(\partial_k \mathbf{v}))\zeta(x/m)|\, dx$$
$$\leq C(1/n) \int_0^T |\nabla \omega_n(s)|_2 |\omega_n(s)|_2\, ds \leq C/n \to 0,$$

as $n \to \infty$, and by (2.23),

$$\int_0^T \int |(\omega_n^k(s)\omega_n(s) - \omega^k(s)\omega(s)|\zeta(x/m)\, dx \to 0.$$

Thus, (2.30) follows. On the other hand,

$$\int (|(\Psi_n(\omega_n^k(s))\omega_n(s))| + |\omega^k(s)\omega(s)||\mathcal{S}(\partial_k \mathbf{v})(1 - \zeta(x/m))|)\, dx$$
$$\leq C(|\omega_n(s)|_2^2 + |\omega(s)|_2^2) \sup_x |\mathcal{S}(\partial_k \mathbf{v})(1 - \zeta(x/m))|$$
$$\leq C \sup_x |\mathcal{S}(\partial_k \mathbf{v})(1 - \zeta(x/m))| \to 0,$$

as $m \to \infty$ (by Sobolev's embedding theorem). So, (2.29) holds and (2.26) is proved. Since the sequence $L_t^{n, \mathbf{v}}(\omega_n)$ is equicontinuous in $t$,

$$\sup_{s \leq T} |L_s^{n, \mathbf{v}}(\omega^n) - L_s^{\mathbf{v}}(\omega)| \to 0. \tag{2.31}$$

Thus for each compact set $K \subseteq \mathcal{Z}$,

$$\sup_{s \leq T, \omega \in K} |L_s^{n, \mathbf{v}}(\omega) - L_s^{\mathbf{v}}(\omega)| \to 0.$$

Since $\{\mathbf{P}^n, n \geq 1\}$ is relatively compact, by Prokhorov's theorem for topological vector spaces (see [1]), for each $\eta > 0$,

$$\lim_n \mathbf{P}^n \left( \sup_{s \leq T} |L_s^{n, \mathbf{v}} - L_s^{\mathbf{v}}| > \eta \right) = 0. \tag{2.32}$$

For $M > 0$, define

$$\tau_M = \inf(t : |L_t^{\mathbf{v}}| > M).$$

Let $\tau^n = \inf(t : |L_t^{n,\mathbf{v}} - L_t^{\mathbf{v}}| > 1)$, $\tau_M^n = \tau_M \wedge \tau^n$. Then we have by (2.32),

(2.33) $$\mathbf{P}^n(\tau^n < t) \leq \mathbf{P}^n\left(\sup_{s \leq T} |L_s^{n,\mathbf{v}} - L_s^{\mathbf{v}}| > 1\right) \to 0,$$

as $n \to \infty$. Also, obviously,

(2.34) $$\sup_{n,s} |L_{s \wedge \tau_L^n}^{n,\mathbf{v}}| \leq M + 1.$$

Let $\mathbf{P}(\tau_M = \tau_{M-}) = 1$, where $\tau_{M-} = \lim_{\delta \downarrow 0} \tau_{M-\delta}$. If $\tau_M(\omega) = \tau_{M-}(\omega)$ and $\omega^n \to \omega$ in $\mathcal{Z}$, we have by (2.31) for each $T > 0$,

(2.35) $$\sup_{s \leq T} |L_{s \wedge \tau_M(\omega^n)}^{\mathbf{v}}(\omega^n) - L_{s \wedge \tau_M(\omega)}^{\mathbf{v}}(\omega)| \to 0.$$

Let $s \leq t$ and let $f$ be a bounded $\mathcal{D}_s$-measurable $\mathbf{P}$-a.s continuous function. Then, by (2.32)–(2.35),

$$0 = \mathbf{P}^n[f(L_{t \wedge \tau_M^n}^{n,\mathbf{v}} - L_{s \wedge \tau_M^n}^{n,\mathbf{v}})] \to \mathbf{P}[f(L_{t \wedge \tau_M}^{\mathbf{v}} - L_{s \wedge \tau_M}^{\mathbf{v}})].$$

Thus, $\mathbf{P}$ is a solution of the martingale problem $(u_0, A, B)$.

In the case $d = 2$, the following inequality holds for all $\mathbf{v} \in \mathbb{H}_2^1$:

(2.36) $$|\mathbf{v}|_4 \leq 2^{1/4} |\mathbf{v}|_2^{1/2} |\nabla \mathbf{v}|_2^{1/2}.$$

By the Hölder inequality, for each $\mathbf{v} \in \mathbb{C}_0^\infty$,

$$\left| \int X^k(t) \, \partial_k X^l(t) v^l \, dx \right| \leq \left( \int |X(t)|^4 \, dx \right)^{1/2} \left( \int |\nabla \mathbf{v}| \, dx \right)^{1/2}$$
$$\leq C |X(t)|_2 |\nabla X(t)|_2 |\mathbf{v}|_{1,2},$$

$dt \times d$, $\mathbf{P}$-a.s. So,

$$|X^k(t) \, \partial_k X(t)|_{-1,2} \leq C |X(t)|_2 |\nabla X(t)|_2,$$

and

$$\int_0^T |X^k(t) \, \partial_k X(t)|_{-1,2}^2 \, dt \leq C \sup_{t \leq T} |X(t)|_2^2 \int_0^T |\nabla X(t)|_2^2 \, dt < \infty,$$

$\mathbf{P}$-a.s. Then, obviously, $\mathbf{P}$-a.s.,

$$\int_0^T |A(t, X_t) - N(t, X_t)|_{-1,2}^2 \, ds < \infty,$$

and by [30], $X_t$ has an $\mathbb{L}_2$-valued (strongly) continuous modification of $X_t$. (Also, the Itô formula holds for $|X_t|_2^2$.) □

Now we shall prove Theorem 2.1.

PROOF OF THEOREM 2.1. According to Theorem 2.10, there is a measure $\mathbf{P}$ on $\mathcal{Z}$ such that (2.21) holds and $\mathbf{P}$-a.s. for each $\mathbf{v} \in \mathbb{C}_0^\infty$,

$$\langle X_t, \mathbf{v} \rangle = \langle \mathbf{u}_0, \mathbf{v} \rangle + \int_0^t \langle A(s, X_s) - N(s, X_s), \mathbf{v} \rangle_{k_0} \, ds + M_t^{\mathbf{v}},$$

where $M_t^{\mathbf{v}} \in \mathcal{M}_{\mathrm{loc}}(\mathbb{D}, \mathbf{P})$ and

$$|M_t^{\mathbf{v}}|^2 - \int_0^t |\langle B(s, X_s), \mathbf{v} \rangle_{k_0, Y}|_Y^2 \, ds \in \mathcal{M}_{\mathrm{loc}}(\mathbb{D}, \mathbf{P}).$$

Since

$$\langle B(s, X_s), \mathbf{v} \rangle_{k_0, Y} = \int (\sigma^i(s) \, \partial_i X_s + \mathbf{g}(s, X_s), \mathcal{S}(\mathbf{v})) \, dx$$

and

$$\mathbf{P} \int_0^T |\sigma^i(s) \, \partial_i X_s + \mathbf{g}(s, X_s)|_2^2 \, ds$$
$$\leq C \mathbf{P} \int_0^T (|\nabla X_s|_2^2 + |X_s|_2^2 + |H(s)|_2^2) \, ds < \infty,$$

there is an $\mathbb{L}_2$-valued continuous martingale $M_t$ such that $\mathbf{P}$-a.s. $\langle M_t, \mathbf{v} \rangle = M_t^{\mathbf{v}}$ for all $t$. Indeed, we simply take an $\mathbb{L}_2$ basis $(\mathbf{e}_k)$ and define

$$M_t = \sum_k M_t^{\mathbf{e}_k} \mathbf{e}_k.$$

Thus, $\mathbf{P}$-a.s.,

(2.37)
$$\partial_t X_t = A(t, X_t) - N(t, X_t) + \partial_t M_t,$$
$$X(0) = \mathbf{u}_0.$$

According to Lemma 3.2 in [25], there exists a cylindrical Wiener process $W$ in $Y$ [possibly in some extension of the probability space $(\Omega, \mathcal{D}_T, \mathbf{P})$] such that

$$M_t = \int_0^t \mathcal{S}(\sigma^i(s, X(s)) \, \partial_i X(s) + \mathbf{g}(t, X(s))) \cdot dW_s,$$

that is, $\partial_t M_t = \mathcal{S}(\sigma^i(s, X(t)) \, \partial_i X(t) + \mathbf{g}(t, X(t))) \cdot \dot{W}_t$. Thus, Theorem 2.1 follows from Theorem 2.10. □

2.4. *Existence and uniqueness of strong global solutions in two dimensions.* To prove Theorem 2.2, we will follow the ideas in [17], where a finite-dimensional stochastic differential equation was considered. First of all, we prove the pathwise uniqueness of the solution in two dimensions.

PROPOSITION 2.11. *Let $d = 2$, $u_0 \in \mathbb{L}_2$, let* (B1), (B2) *hold and for all $l, j, t, x$ and every $\mathbf{u}, \bar{\mathbf{u}}$,*

$$|f^l(t, x, \mathbf{u}) - f^l(t, x, \bar{\mathbf{u}})| + |g^l(t, x, \mathbf{u}) - g^l(t, x, \bar{\mathbf{u}})|_Y$$
$$+ |f^{l,j}(t, x, \mathbf{u}) - f^{l,j}(t, x, \bar{\mathbf{u}})|$$
$$\leq K|\mathbf{u} - \bar{\mathbf{u}}|.$$

*Assume that on some probability space $(\Omega, \mathcal{F}, \mathbf{P})$, with a right-continuous filtration of $\sigma$-algebras $\mathbb{F} = (\mathcal{F}_t)$ and cylindrical Wiener process $W$ in $Y$, we have two solutions $\mathbf{U}_1, \mathbf{U}_2$ to the Navier–Stokes equation* (2.3) *such that $\mathbf{P}$-a.s.,*

$$\sup_{s \leq T} |\mathbf{U}_l(s)|_2^2 + \int_0^T |\nabla \mathbf{U}_l(s)|_2^2 \, ds < \infty, \qquad l = 1, 2.$$

*Then $\mathbf{P}$-a.s., $\mathbf{U}_1(t) = \mathbf{U}_2(t)$ for all $t$.*

PROOF. Let $\mathbf{U} = \mathbf{U}_1 - \mathbf{U}_2$. We apply the Itô formula for $|\mathbf{U}(t)|_2^2$:

(2.38)
$$|\mathbf{U}(t)|_2^2 = |\mathbf{U}(0)|_2^2 + \int_0^t [2\langle A(s, \mathbf{U}_2(s)) - A(s, \mathbf{U}_1(s)), \mathbf{U}(s)\rangle_{k_0}$$
$$- 2\langle N(s, \mathbf{U}_2(s)) - N(s, \mathbf{U}_1(s)), \mathbf{U}(s)\rangle_{k_0}$$
$$+ \|B(s, \mathbf{U}_2(s)) - B(s, \mathbf{U}_1(s))\|_2^2] \, ds$$
$$+ 2 \int_0^t \langle B(s, \mathbf{U}_2(s)) - B(s, \mathbf{U}_1(s)), \mathbf{U}(s)\rangle_{0,Y} \cdot dW_s.$$

Notice

$$\int (U_1^k(t) \, \partial_k U_1^l(t) - U_2^k(t) \, \partial_k U_2^l(t))(U_1^l(t) - U_2^l(t)) \, dx$$
$$= \int (U_1^k(t) - U_2^k(t)) \partial_k U_2^l(t)(U_1^l(t) - U_2^l(t)) \, dx$$
$$= \int U^k(t) \, \partial_k U_2^l(t) U^l(t) \, dx.$$

Using (2.36), we find that for each $\varepsilon$ there is a constant $C_\varepsilon$ such that

(2.39)
$$\left| \int (U_1^k(t) \, \partial_k U_1^l(t) - U_2^k(t)\partial_k U_2^l(t))(U_1^l(t) - U_2^l(t)) \right| dx$$
$$\leq \left( \int |\mathbf{U}(t)|^4 \, dx \right)^{1/2} \left( \int |\nabla \mathbf{U}_2(t)|^2 \, dx \right)^{1/2}$$
$$\leq C(|\mathbf{U}(t)|_2 |\nabla \mathbf{U}(t)|_2 |\nabla \mathbf{U}_2|_2) \leq \varepsilon |\nabla \mathbf{U}|_2^2 + C_\varepsilon |\mathbf{U}(t)|_2^2 |\nabla \mathbf{U}_2(t)|_2^2.$$

Let $\tau$ be an arbitrary stopping time such that

$$\mathbf{E}\int_0^\tau \left| \int [\sigma^{l,j}(s)\partial_j U^l(s) + g^l(s, \mathbf{U}_2(s)) - g^l(s, \mathbf{U}_2(s))]U^l(s)\, dx \right|_Y^2 ds < \infty.$$

By (2.38), (2.39) and our assumptions, there are some constants $\varepsilon$ and $C$ independent of $\tau$ such that for all $t$,

$$\mathbf{E}\left[|\mathbf{U}(t\wedge\tau)|_2^2 + \varepsilon\int_0^{\tau\wedge t}|\nabla \mathbf{U}(s)|_2^2\, ds\right] \leq C\mathbf{E}\int_0^{t\wedge\tau}|\mathbf{U}(s)|_2^2\, d\lambda_s,$$

where $\lambda_s = s + \int_0^s |\nabla \mathbf{U}_2(r)|_2^2\, dr$. The pathwise uniqueness now follows (see, e.g., Lemma 2 in [16]). □

Let $(\Omega, \mathcal{F}, \mathbf{P})$ be a probability space with a right-continuous filtration of $\sigma$-algebras $\mathbb{F} = (\mathcal{F}_t)$ and a cylindrical Wiener process $W$ in $Y$. Let (B1) and (B2) be satisfied and let $E|u_0|_2^2 < \infty$. Then, according to Proposition 2.4, for each $n$, there exists a unique $\mathbb{L}_2$-valued continuous solution $u_n(t)$ of (2.7) such that $\int_0^T |\nabla u_n(s)|_2^2\, ds < \infty$ and $\operatorname{div} u_n(t) = 0$ for all $t$, $\mathbf{P}$-a.s. Moreover,

$$\sup_n \mathbf{E}\left[\sup_{t\leq T}|\mathbf{u}_n(t)|_2^2 + \int_0^T |\nabla \mathbf{u}_n(t)|_2^2\, dt\right] < \infty.$$

PROPOSITION 2.12. *Assume $d = 2$, $u_0 \in \mathbb{L}_2$, (B1), (B2) hold and for all $l, j, t, x$ and every $\mathbf{u}, \bar{\mathbf{u}}$,*

$$|f^l(t, x, \mathbf{u}) - f^l(t, x, \bar{\mathbf{u}})| + |g^l(t, x, \mathbf{u}) - g^l t, x, \bar{\mathbf{u}})|_Y$$
$$+ |f^{l,j}(t, x, \mathbf{u}) - f^{l,j}(t, x, \bar{\mathbf{u}})|$$
$$\leq K|\mathbf{u} - \bar{\mathbf{u}}|.$$

*Then there exists a unique $\mathbb{L}_2$-valued continuous solution $u(t)$, $t \in [0, T]$, of (2.3) on $(\Omega, \mathcal{F}, \mathbf{P})$ such that $\int_0^T |\nabla u(s)|_2^2\, ds < \infty$, $\operatorname{div} u(t) = 0$ for all $t$ $\mathbf{P}$-a.s. Moreover, for each $\mathbf{v} \in \mathbb{C}_0^\infty(\mathbf{R}^2)$, $i = 1, 2$, $R > 0$,*

$$\sup_t |\langle \mathbf{u}_n(t) - \mathbf{u}(t), \mathbf{v}\rangle_0| + \int_0^T \int_{|x|\leq R}|\mathbf{u}_n(t,x) - \mathbf{u}(t,x)|^2\, dx\, dt$$
(2.40)
$$+ \left|\int_0^T \langle \partial_i \mathbf{u}_n(t) - \partial_i \mathbf{u}(t), \mathbf{v}\rangle_0\, dt\right| \to 0$$

*in probability, as $n \to \infty$.*

PROOF. Let $(e_n)$ be a CONS of the separable Hilbert space $Y$. Then,

$$\tilde{Y} = \left\{ y = \sum_{k=1}^\infty y_k e_k : y_k \in \mathbf{R}, \sum_{k=1}^\infty \frac{y_k^2}{k^2} < \infty \right\}$$

is a Hilbert space with the norm

$$|y|_{\tilde{Y}} = \left(\sum_k \frac{y_k^2}{k^2}\right)^{1/2},$$

and $W_t$ is a $\tilde{Y}$-valued continuous process. Consider

$$E = \mathcal{Z} \times \mathcal{Z} \times C_{[0,T]}(\tilde{Y})$$

with the product of corresponding topologies and denote $(X_t^1, X_t^2, \overline{W}_t)$ the canonical process in $E$. For $\tilde{\omega} = (\omega, \bar{\omega}, w) \in E$,

$$X_t^1 = X^1(t, \tilde{\omega}) = \omega(t), \qquad X_t^2 = X^2(t, \tilde{\omega}) = \bar{\omega}(t), \qquad \overline{W}_t(\tilde{\omega}) = w(t).$$

Let $\mathcal{E}_t = \sigma(X_s^1, s \leq t) \otimes \sigma(X_s^2, s \leq t) \otimes \sigma(\overline{W}_s, s \leq t)$, $\mathbb{E} = (\mathcal{E}_{t+})_{0 \leq t < T}$, $\mathcal{E} = \mathcal{E}_T$. For each $m, l$, the process $(u_m, \mathbf{u}_l, W)$ induces a measure $\mathbf{P}^{m,l}$ on $(E, \mathcal{E})$. For $\mathbf{v} \in \mathbb{C}_0^\infty(\mathbf{R}^d)$, $k = 1, 2$, $y \in Y$, denote

$$\varphi_{k,n}^{\mathbf{v},y}(s, X^k) = i_0 \langle A_n(s, X_s^k) - N_n(s, X_s^k), \mathbf{v}\rangle_{k_0} - \tfrac{1}{2}|\langle B_n(s, X_s^k), \mathbf{v}\rangle_{k_0, Y} + y|_Y^2,$$

where $i_0^2 = -1$. Let

$$L_t^{n,k,\mathbf{v},y} = e^{i_0(\langle X_t^k, \mathbf{v}\rangle + \overline{W}_t(y))} - \int_0^t e^{i_0(\langle X_s^k, \mathbf{v}\rangle + \overline{W}_s(y))} \varphi_{k,n}^{\mathbf{v},y}(s, X^k)\, ds,$$

$k = 1, 2$, $n \geq 1$. Applying the Itô formula to the scalar semimartingale $\langle X_t^k, \mathbf{v}\rangle + W_t(y) = \langle X_t^k, \mathbf{v}\rangle_0 + W_t(y)$, $k = 1, 2$, we obtain that for each test function $\mathbf{v} \in \mathbb{C}_0^\infty(\mathbf{R}^d)$, $y \in Y$,

(2.41) $$L_t^{m,1,\mathbf{v},y}, L_t^{l,2,\mathbf{v},y} \in \mathcal{M}_{\mathrm{loc}}^c(\mathbb{E}, \mathbf{P}^{m,l}).$$

Since the set $\{\mathbf{P}^n, n \geq 1\}$ of probability measures on $\mathcal{Z}$ is relatively compact (see Corollary 2.8), the set $\{\mathbf{P}^{m,l} : m \geq 1, l \geq 1\}$ is relatively compact. Assume that for some subsequences $m(n) \to \infty$, $l(n) \to \infty$, $\mathbf{P}^{m(n), l(n)} \to \overline{\mathbf{P}}$, as $n \to \infty$. We will prove now that

(2.42) $$\overline{\mathbf{P}}(X_t^1 = X_t^2, t \in [0, T]) = 1.$$

Obviously, $\overline{W}_t$ is a cylindrical Wiener process in $Y$ with respect to the filtration $\mathbb{E}$ on the probability space $(E, \mathcal{E}, \overline{\mathbf{P}})$. Let

$$\varphi_k^{\mathbf{v},y}(s, X_s^k) = i_0 \langle A(s, X_s^k) - N(s, X_s^k), \mathbf{v}\rangle_{k_0} - \tfrac{1}{2}|\langle B(s, X_s^k), \mathbf{v}\rangle_{k_0, Y} + y|_Y^2,$$

$k = 1, 2$, $\mathbf{v} \in \mathbb{C}_0^\infty(\mathbf{R}^d)$, $y \in Y$. Define

$$L_t^{k,\mathbf{v},y} = e^{i_0(\langle X_t, \mathbf{v}\rangle + \overline{W}_t(y))} - \int_0^t e^{i_0(\langle X(s), \mathbf{v}\rangle + \overline{W}_s(y))} \varphi_k^{\mathbf{v},y}(s, X_s)\, ds.$$

In a standard way, (2.41) implies (see Theorem 2.10) that

(2.43) $$L_t^{1,\mathbf{v},y}, L_t^{2,\mathbf{v},y} \in \mathcal{M}_{\mathrm{loc}}^c(\mathbb{E}, \overline{\mathbf{P}}),$$

which means that both $X_t^1$ and $X_t^2$ satisfy (2.3). Indeed, (2.43) implies that for all $\mathbf{v} \in \mathbb{C}_0^\infty(\mathbf{R}^d)$,

$$\langle X_t^k, \mathbf{v} \rangle = \langle \mathbf{u}_0, \mathbf{v} \rangle + \int_0^t \langle A(s, X_s^k) - N(s, X_s^k), \mathbf{v} \rangle_{k_0} \, ds + M_t^{k,\mathbf{v}},$$

where $M_t^{k,\mathbf{v}} \in \mathcal{M}_{\mathrm{loc}}(\mathbb{E}, \overline{\mathbf{P}})$,

$$|M_t^{k,\mathbf{v}}|^2 - \int_0^t |\langle B(s, X_s^k), \mathbf{v} \rangle_{k_0, Y}|_Y^2 \, ds \in \mathcal{M}_{\mathrm{loc}}(\mathbb{E}, \overline{\mathbf{P}}),$$

and finally

$$M_t^{k,\mathbf{v}} = \int_0^t \langle B(s, X_s^k), \mathbf{v} \rangle_{k_0, Y} \, d\overline{W}_s,$$

$k = 1, 2$. Therefore, by Proposition 2.11, (2.42) holds.

Let

$$I_{m,l}^1 = \int_0^T \int_{|x| \leq R} |\mathbf{u}_m(t,x) - \mathbf{u}_l(t,x)|^2 \, dx \, dt,$$

$$I_{m,l}^2 = \sup_t |\langle \mathbf{u}_m(t) - \mathbf{u}_l(t), \mathbf{v} \rangle|,$$

$$I_{m,l}^3 = \left| \int_0^T \langle \partial_i \mathbf{u}_m(t) - \partial_i \mathbf{u}_l(t), \mathbf{v} \rangle_0 \, dt \right|,$$

$R > 0$, $\mathbf{v} \in \mathbb{C}_0^\infty(\mathbf{R}^d)$, $i = 1, 2$. Let $F(x) = |x| \wedge 1$, $x \in \mathbf{R}$. From the weak convergence of $\mathbf{P}^{m(n),l(n)}$ to $\overline{\mathbf{P}}$ and (2.42), it follows that

$$\mathbf{E} F(I_{m(n),l(n)}^j) \to 0,$$

as $n \to \infty$, $j = 1, 2, 3$. Since this is true for an arbitrary converging subsequence, we have the convergence in probability of $I_{m,l}^j$, $j = 1, 2, 3$, as $m, l \to \infty$. Therefore, there exists an $\mathbb{L}_2$-valued weakly continuous process $u(t)$, $t \in [0, T]$, on $(\Omega, \mathcal{F}, \mathbf{P})$ such that $\int_0^T |\nabla u(s)|_2^2 \, ds < \infty$, $\mathrm{div}\, u(t) = 0$ for all $t$, $\mathbf{P}$-a.s. and (2.40) holds. For each $n$ and $\mathbf{v} \in \mathbb{C}_0^\infty(\mathbf{R}^d)$,

$$\partial_t \langle \mathbf{u}_n(t), \mathbf{v} \rangle = \langle A_n(t, \mathbf{u}_n(t)), \mathbf{v} \rangle - \langle N_n(t, \mathbf{u}_n(t)), \mathbf{v} \rangle + \langle B_n(t, \mathbf{u}_n(t)), \mathbf{v} \rangle_{0,Y} \cdot \dot{W}_t,$$

$$\mathbf{u}_n(0) = \mathbf{u}_{0,n}, \qquad \mathrm{div}\, \mathbf{u}_n(t) = 0,$$

and passing to the limit, as $n \to \infty$, in this equation we see that $\mathbf{u}(t)$ satisfies (2.3). According to Theorem 2.10, $\mathbf{u}(t)$ is strongly continuous in $\mathbb{L}_2$ and the statement follows. $\square$

Now we can prove Theorem 2.2.

PROOF OF THEOREM 2.2. Assume that on some probability space $(\Omega, \mathcal{F}, \mathbf{P})$, with a right-continuous filtration of $\sigma$-algebras $\mathbb{F} = (\mathcal{F}_t)$ and cylindrical

Wiener process $W$ in $Y$, we have a solution $\mathbf{u}$ to the Navier–Stokes equation (2.3) such that $\mathbf{P}$-a.s.,

$$\sup_{s\leq T}|\mathbf{u}(s)|_2^2 + \int_0^T |\nabla \mathbf{u}(s)|_2^2\, ds < \infty.$$

By Propositions 2.11 and 2.12, we have the convergence (2.40) for the approximating sequence $\mathbf{u}_n(t)$. Since $\mathbf{u}_n(t)$ can be constructed by iterations (see [23]), the distribution of all $\mathbf{u}_n(t)$ and therefore the distribution of $\mathbf{u}(t)$ are uniquely determined by $\mathbf{u}_0$. $\square$

**3. Wiener chaos and strong solutions.** In this section we will derive a system of deterministic PDEs for Fourier coefficients of the Wiener chaos expansion of a solution of stochastic Navier–Stokes equation (2.1). This system is usually referred to as the propagator. We will demonstrate that the existence of a solution of the propagator is a necessary and sufficient condition for the existence of a strong (pathwise) solution of the related Navier–Stokes equation. Similarly, the uniqueness of a solution of the propagator is equivalent to the pathwise uniqueness of the related stochastic Navier–Stokes equation.

First we shall introduce additional notation and recall some basic facts of the Wiener chaos theory (see, e.g., [18, 19, 22], etc.).

Let us fix a positive number $T < \infty$. Let $\{m_k, k \geq 1\}$ be an orthonormal basis in $L_2(0,T)$ and let $\{\ell_k, k \geq 1\}$ be an orthonormal basis in $Y$. Write $\xi_i^k = \int_0^T m_i(s)\, dw^k(t)$, where $w_k(t) = (W(t), \ell_k)_Y$. Let $\alpha = \{\alpha_i^k, k \geq 0; i \geq 1\}$ be a multi-index; that is, for every $(i,k)$, $\alpha_i^k \in \mathbb{N} = \{0, 1, 2, \dots\}$. We shall consider only such $\alpha$ that $|\alpha| = \sum_{k,i} \alpha_i^k < \infty$, that is, only a finite number of $\alpha_i^k$ are nonzero, and we denote by $\mathcal{J}$ the set of all such multi-indices. Obviously, if $\alpha \in \mathcal{J}$, the number $\alpha! = \prod_{k,i} \alpha_i^k!$ is well defined. For $\alpha \in \mathcal{J}$, write

$$\zeta_\alpha := \prod_{i,k=1}^\infty H_{\alpha_i^k}(\xi_i^k),$$

where $H_n$ is the $n$th Hermite polynomial. The random variable $\zeta_\alpha$ is often referred to as the (unnormalized) $\alpha$th Wick polynomial.

The most important feature of the Wick polynomials $\zeta_\alpha$ is that the set $\{\zeta_\alpha/\sqrt{\alpha!}, \alpha \in \mathcal{J}\}$ is an orthonormal basis in $L_2(\Omega, \mathcal{F}_T, \mathbf{P})$, where $\mathcal{F}_t = \sigma(W(s),\, s \leq t)$ (see, e.g., [4, 22]). This result is often referred to as the Cameron–Martin theorem.

By Lemma 15 in [29], the process $\zeta_\alpha(t) = E[\zeta_\alpha | \mathcal{F}_t]$ satisfies the following equation:

$$(3.1) \qquad d\zeta_\alpha(t) = D\zeta_\alpha(t) \cdot dW(t),$$

where $D\zeta_\alpha(t) = m_i(t)\alpha_i^k \zeta_{\alpha(i,k)}(t)\ell_k$ is the Malliavin derivative of $\zeta_\alpha(t)$; the multi-index $\alpha(i,j) \in \mathcal{J}$ is defined by

$$\alpha(i,j)_l^k = \begin{cases} \alpha_l^k, & \text{if } (k,l) \neq (j,i) \text{ or } k = 0, \\ (\alpha_l^k - 1) \vee 0, & \text{if } (k,l) = (j,i). \end{cases} \tag{3.2}$$

For $\alpha, \beta \in \mathcal{J}$, define $|\alpha - \beta| = (|a_1 - \beta_1|, |a_2 - \beta_2|, \dots)$.

DEFINITION 3.1 (cf. [29]). We say that a triple of multi-indices $(\alpha, \beta, \gamma)$ is complete, written $(\alpha, \beta, \gamma) \in \mathcal{C}$, if all the entries of the multi-index $\alpha + \beta + \gamma$ are even numbers and $|\alpha - \beta| \leq \gamma \leq \alpha + \beta$.

It is readily checked that the following criterion holds:

LEMMA 3.1. *A triple $(\alpha, \beta, \gamma)$ is complete if and only if $\alpha + \beta + \gamma = 2p$ for some $p \in \mathcal{J}$ and $p \leq \alpha \wedge \beta$.*

For $(\alpha, \beta, \gamma) \in \mathcal{C}$, we define

$$\Phi(\alpha, \beta, \gamma) = \left(\left(\frac{\alpha - \beta + \gamma}{2}\right)! \left(\frac{\beta - \alpha + \gamma}{2}\right)! \left(\frac{\alpha + \beta - \gamma}{2}\right)!\right)^{-1}.$$

Obviously, $\Phi(\alpha, \beta, \gamma)$ is invariant with respect to permutations of the arguments.

For $\alpha \in \mathcal{J}$, write $U^\alpha = \{\gamma, \beta \in \mathcal{J} : (\alpha, \beta, \gamma) \in \mathcal{C}\}$.

Now we can derive the Wiener chaos expansion for a strong solution of stochastic Navier–Stokes equation (2.1).

For the sake of simplicity, in addition to assumptions of Section 2.1, throughout this section we will assume

(C1) Functions $f^l = f^l(t,x)$, $f^{l,j} = f^{l,j}(t,x)$ and $g^l = g^l(t,x)$ do not depend on $\mathbf{u}$.

THEOREM 3.2. *Let $d \geq 2$ and $u_0 \in \mathbb{L}_2$. Assume that (2.1) has a weakly continuous strong solution $u(t)$ such that*

$$\sup_{s \leq T} \mathbf{E}|\mathbf{u}(s)|_2^2 + \int_0^T \mathbf{E}|\nabla \mathbf{u}(s)|_2^2 \, ds < \infty. \tag{3.3}$$

*Then*

$$\mathbf{u}(t) = \sum_{\alpha \in \mathcal{J}} \frac{\hat{\mathbf{u}}_\alpha(t)}{\sqrt{\alpha!}} \xi_\alpha, \tag{3.4}$$

*and the Hermite–Fourier coefficients $\hat{\mathbf{u}}_\alpha(t)$ are $\mathbb{L}_2$-valued weakly continuous functions so that*

$$\sup_{s \leq T} \sum_{\alpha \in \mathcal{J}} \frac{|\hat{\mathbf{u}}_\alpha(s)|_2^2}{\alpha!} + \int_0^T \sum_{\alpha \in \mathcal{J}} \frac{|\nabla \hat{\mathbf{u}}_\alpha(s)|_2^2}{\alpha!} \, ds < \infty. \tag{3.5}$$

*Moreover, the set of functions $\{\hat{\mathbf{u}}_\alpha(t,x), \alpha \in \mathcal{J}\}$ satisfies the propagator equation*

$$
\begin{aligned}
\partial_t \hat{\mathbf{u}}_\alpha(t) = \mathcal{S}\Bigg[ & \partial_i(a^{ij}(t)\,\partial_j \hat{\mathbf{u}}_\alpha(t)) \\
& - \sum_{\gamma,\beta \in U^\alpha} \alpha!\hat{u}^i_\gamma(t)\,\partial_i \hat{\mathbf{u}}_\beta(t)\Phi(\alpha,\beta,\gamma) \\
& + b^i(t)\,\partial_i \hat{\mathbf{u}}_\alpha(t) + \tilde{L}_i(\hat{\mathbf{u}}_\alpha(t),t)\mathbf{h}^i(t) \\
& + I_{\{|\alpha|=0\}}(\mathbf{f}(t) + \partial_j \mathbf{f}^j(t, \mathbf{u}(t))) \\
& + m_j(t)\alpha_j^k(\sigma^{ik}(t)\,\partial_i \hat{\mathbf{u}}_{\alpha(j,k)}(t) + I_{\{|\alpha|=1\}}\mathbf{g}_k(t))\Bigg], \qquad t\in(0,T];
\end{aligned}
\tag{3.6}
$$

$$\hat{\mathbf{u}}_\alpha(0) = I_{\{|\alpha|=0\}}\mathbf{u}_0,$$

*where* $\mathbf{g}_k(t) = (\mathbf{g}(t), \ell_k)_Y$, *and*

$$\tilde{\mathbf{L}}(\mathbf{u}(t),t) = (\tilde{L}_l(\mathbf{u}(t),t))_{1\leq l \leq d} = \mathcal{G}(\sigma^k(t)\,\partial_k\mathbf{u}(t) + \mathbf{g}(t)).$$

PROOF. Let us fix a complete orthonormal system $\{\mathbf{e}_k, k\geq 1\}$ in $\mathbb{L}_2$ so that every $\mathbf{e}_k \in \mathbb{H}_2^1$.

Let $u$ be a strong global solution of (2.1) which is $\mathbb{L}_2$-weakly continuous and so that (3.3) holds. By the Cameron–Martin theorem, for every $t$,

$$\hat{\mathbf{u}}_\alpha(t) = \sum_{i=1}^\infty \mathbf{E}[(\mathbf{u}(t), \mathbf{e}_i)_2 \zeta_\alpha(t)]\mathbf{e}_i \tag{3.7}$$

and

$$
\begin{aligned}
\mathbf{E}|\mathbf{u}(t)|_2^2 &= \sum_{\alpha \in \mathcal{J}} \frac{1}{\alpha!}\sum_{i=1}^\infty (\mathbf{E}[(\mathbf{u}(t),\mathbf{e}_i)_2\zeta_\alpha(t)])^2 \\
&= \sum_{\alpha \in \mathcal{J}} \frac{1}{\alpha!}|\hat{\mathbf{u}}_\alpha(t)|_2^2 < \infty.
\end{aligned}
$$

Let us fix $\alpha \in \mathcal{J}$. Owing to (3.3), for any $\mathbf{v} \in \mathbb{L}_2$ and any set $\{t_n\}$ of points in $[0,T]$, the sequence $\{\zeta_\alpha(\mathbf{u}(t_n), \mathbf{v})\}$ is uniformly $P$-integrable. Therefore, the weak continuity of $u(t)$ implies that $\hat{\mathbf{u}}_\alpha(t)$ is also weakly continuous in $\mathbb{L}_2$.

The same arguments as before yield that for almost all $t$,

$$\widehat{(\partial_j \mathbf{u}(t))}_\alpha := \sum_{i=1}^\infty \mathbf{E}[(\partial_j\mathbf{u}(t), \mathbf{e}_i)\zeta_\alpha(t)]\mathbf{e}_i \tag{3.8}$$

and

$$\int_0^T \sum_{\alpha \in \mathcal{J}} \frac{1}{\alpha!}|\widehat{(\partial_j\mathbf{u})}_\alpha(t)|_2^2\,dt = \int_0^T \mathbf{E}|\partial_j\mathbf{u}(t)|_2^2\,dt < \infty.$$

By (3.3) and the Fubini theorem, for all $\varphi \in \mathbb{C}_0^\infty(\mathbf{R}^d)$ and almost all $t \leq T$,

$$\int_{\mathbf{R}^d} ((\widehat{\partial_i \mathbf{u}(t)})_\alpha, \varphi) \, dx = -\int_{\mathbf{R}^d} \mathbf{E}(\mathbf{u}(t), \partial_i \varphi) \zeta_\alpha \, dx = -\int_{\mathbf{R}^d} (\hat{\mathbf{u}}_\alpha(t), \partial_i \varphi) \, dx.$$

Thus, for almost all $t \leq T$,

(3.9) $$\partial_i \hat{\mathbf{u}}_\alpha(t) = (\widehat{\partial_i \mathbf{u}(t)})_\alpha$$

and

(3.10) $$\int_0^T \sum_{\alpha \in \mathcal{J}} \frac{1}{\alpha!} |\hat{\mathbf{u}}_\alpha(t)|_{1,2}^2 \, dt = \int_0^T \mathbf{E}|\mathbf{u}(t)|_{1,2}^2 \, dt < \infty.$$

Denote $\mathcal{V} := \{\varphi : \varphi \in \mathbb{C}_0^\infty(\mathbf{R}^d), \operatorname{div} \varphi = 0\}$. Obviously, for every $t$ in $(0, T]$ and every test function $\varphi \in \mathcal{V}$, $P$-a.a.,

$$\int_{\mathbf{R}^d} (\mathbf{u}(t), \varphi) \, dx$$

$$= \int_{\mathbf{R}^d} (\mathbf{u}_0, \varphi) \, dx - \int_0^t \int_{\mathbf{R}^d} a^{ij}(s)(\partial_j \mathbf{u}(s), \partial_i \varphi) \, dx \, ds$$

(3.11) $$+ \int_0^t \int_{\mathbf{R}^d} (\varphi, -u^k(s) \, \partial_k \mathbf{u}(s) + b^i(s) \, \partial_i \mathbf{u}(s) + \tilde{L}_i(\mathbf{u}(s), s) \mathbf{h}^i(s)$$

$$+ \mathbf{f}(s) + \partial_i \mathbf{f}^i(s)) \, dx \, ds$$

$$+ \int_0^t \int_{\mathbf{R}^d} (\varphi, \sigma^i(s) \partial_i \mathbf{u}(s) + \mathbf{g}(s)) \, dx \, dW_s.$$

By (3.1) and the Itô formula, we have

$$d\zeta_\alpha(t) \int_{\mathbf{R}^d} (\mathbf{u}(t), \varphi) \, dx$$

$$= \int_{\mathbf{R}^d} (a^{ij}(t) \, \partial_j \mathbf{u}(t) \zeta_\alpha(t), \partial_i \varphi) \, dx \, dt$$

$$+ \int_{\mathbf{R}^d} \zeta_\alpha(t)(\varphi, -u^k(t) \, \partial_k \mathbf{u}(t) + b^i(t) \, \partial_i \mathbf{u}(t)$$

$$+ \tilde{L}_i(\mathbf{u}(t), t) \mathbf{h}^i(t) + \mathbf{f}(t) + \partial_i \mathbf{f}^i(t)) \, dx \, dt$$

$$+ \zeta_\alpha(t) \int_{\mathbf{R}^d} (\varphi, \sigma^i(t) \, \partial_i \mathbf{u}(t) + \mathbf{g}(t)) \, dx \, dW_t$$

$$+ \int_{\mathbf{R}^d} (\varphi, \mathbf{u}(t) D \zeta_\alpha(t)) \, dx \, dW_t$$

$$+ \int_{\mathbf{R}^d} (m_j(t) \alpha_j^k \varphi, (\sigma^{ik}(t) \, \partial_i \hat{\mathbf{u}}_{\alpha(j,k)}(t) + I_{\{|\alpha|=1\}} \mathbf{g}_k(t))) \, dx \, dt,$$

where $\sigma^{ik} = (\sigma^i, \ell_k)_Y$.

This together with (3.9) and (3.10) yields

$$
\begin{aligned}
&\int_{\mathbf{R}^d} (\hat{\mathbf{u}}_\alpha(t), \varphi) \, dx \\
&= I_{\{|\alpha|=0\}} \int_{\mathbf{R}^d} (\mathbf{u}(0), \varphi) \, dx \\
&\quad + \int_0^t \int_{\mathbf{R}^d} \{ (a^{ij}(s)\, \partial_j \hat{\mathbf{u}}_\alpha(s), \partial_i \varphi) \\
&\qquad + (\varphi, -(u^k \widehat{(s)\, \partial_k} \mathbf{u}(s))_\alpha \\
&\qquad + b^i(s)\, \partial_i \hat{\mathbf{u}}_\alpha(s) + \tilde{L}_i(\hat{\mathbf{u}}_\alpha(s), t)\mathbf{h}^i(s) \\
&\qquad + I_{\{|\alpha|=0\}}(\mathbf{f}(s) + \partial_i \mathbf{f}^i) \\
&\qquad + m_j(s)\alpha_j^k (\sigma^{ik}(t)\, \partial_i \hat{\mathbf{u}}(s)_{\alpha(j,k)} + I_{\{|\alpha|=1\}}\mathbf{g}(s))) \} \, dx \, ds.
\end{aligned}
$$
(3.12)

Mimicking the derivation in [29] [see (4.20) and (4.21)], one can easily demonstrate that $ds\, dx$-a.s.,

$$
(u^k \widehat{(s)\, \partial_k} \mathbf{u}(s))_\alpha = \sum_{\gamma, \beta \in U^\alpha} \hat{u}^i_\gamma(s)\, \partial_i \hat{\mathbf{u}}_\beta(s) \Phi(\gamma, \beta, \alpha) \alpha!.
$$
(3.13)

This completes the proof. □

Obviously, Theorem 3.2 and the Cameron–Martin theorem yield the following result:

COROLLARY 3.3. *The first and second moments of a solution of* (2.1) *are given by* $Eu(t) = \hat{\mathbf{u}}_0(t)$ *and*

$$
\mathbf{E}|\mathbf{u}(t)|^2 = \sum_{\alpha \in \mathcal{J}} \frac{|\hat{\mathbf{u}}_\alpha(t)|^2}{\alpha!}.
$$

Now we will prove that the existence of a solution of the propagator equation is not only necessary but also sufficient for the existence of a strong solution of a turbulent stochastic Navier–Stokes equation.

THEOREM 3.4. *Let $d \geq 2$ and $u_0 \in \mathbb{L}_2$. Assume that equation* (3.6) *has a solution $\{\hat{\mathbf{u}}_\alpha(t,x),\ \alpha \in \mathcal{J}\}$ on the interval $(0,T]$ so that for every $\alpha$, $\hat{\mathbf{u}}_\alpha(t)$ is weakly continuous in $\mathbb{L}_2$ and the inequality*

$$
\sup_{s \leq T} \sum_{\alpha \in \mathcal{J}} \frac{|\hat{\mathbf{u}}_\alpha(s)|_2^2}{\alpha!} + \int_0^T \sum_{\alpha \in \mathcal{J}} \frac{|\nabla \hat{\mathbf{u}}_\alpha(s)|_2^2}{a!} \, ds < \infty
$$

*holds.*

Then any $\mathcal{F}_t$-adapted process of the form

(3.14) $$\bar{\mathbf{u}}(t) := \sum_{\alpha \in \mathcal{J}} \frac{\hat{\mathbf{u}}_\alpha(t)}{a!} \zeta_\alpha$$

is a global strong solution of (2.1). Moreover, for every $\mathbf{v} \in \mathbb{L}_2$, $(\bar{\mathbf{u}}(t), \mathbf{v})_2$ is continuous in $L_2(\Omega)$ and

$$\sup_{s \leq T} \mathbf{E} |\bar{\mathbf{u}}(s)|_2^2 + \int_0^T \mathbf{E} |(\nabla \bar{\mathbf{u}}(s))|_2^2 \, ds < \infty.$$

PROOF. By (3.5), $\sup_{t \leq T} \mathbf{E} |\bar{\mathbf{u}}(t)|_2^2 < \infty$. Next, we shall prove that $\mathbf{E} \times \int_0^T |\nabla \bar{\mathbf{u}}(t)|_2^2 \, dt < \infty$.

Write

$$\|\mathbf{v}\|_{1,2} := \left( \int_{\mathbf{R}^d} (1 + |\lambda|^2) |\mathcal{F}\mathbf{v}(t, \lambda)|^2 \, d\lambda \right)^{1/2},$$

where $\mathcal{F}v(\lambda) := \frac{1}{(2\pi)^{d/2}} \int_{\mathbf{R}^d} \exp\{-i(x, \lambda)\} v(x) \, dx$.

Since the norms $\|\cdot\|_{1,2}$ and $|\mathbf{v}|_{1,2}$ are equivalent, we have

$$\mathbf{E} \int_0^T |(\nabla \bar{\mathbf{u}}(t))|_2^2 \, dt \leq CE \int_0^T \|\bar{\mathbf{u}}(t)\|_{1,2}^2 \, dt$$

$$= C \int_0^T \sum_{\alpha \in \mathcal{J}} \|\hat{\mathbf{u}}_\alpha(t)\|_{1,2}^2 / \alpha! \, dt$$

$$\leq C' \int_0^T \sum_{\alpha \in \mathcal{J}} |\hat{\mathbf{u}}_\alpha(t)|_{1,2}^2 / \alpha! \, dt < \infty.$$

Now let us show that, for every $\mathbf{v} \in \mathbb{L}_2$, $(\bar{\mathbf{u}}(t), \mathbf{v})_{\mathbb{L}_2}$ is a mean-square continuous stochastic process. Indeed, let $\{\mathcal{J}_N\}$ be an increasing system of finite subsets of $\mathcal{J}$ so that $\mathcal{J}_N \uparrow \mathcal{J}$. Write $\bar{\mathbf{u}}^N(t) := \sum_{\alpha \in \mathcal{J}_N} \hat{\mathbf{u}}_\alpha(t) \xi_\alpha / \sqrt{\alpha!}$. For every $N$, $\sum_{\alpha \in \mathcal{J}_N} (\hat{\mathbf{u}}_\alpha(t), \mathbf{v})_{\mathbb{L}_2}^2$ is a continuous function of $t$; consequently, $(\bar{\mathbf{u}}^N(t), \mathbf{v})_{\mathbb{L}_2}$ is continuous in $L_2(\Omega)$. Moreover,

$$\lim_{N \to \infty} \sup_t \mathbf{E} |(\bar{\mathbf{u}}^N(t), \mathbf{v})_{\mathbb{L}_2} - (\bar{\mathbf{u}}(t), \mathbf{v})_{\mathbb{L}_2}|^2$$

$$= \lim_{N \to \infty} \sup_t \sum_{\alpha \in \mathcal{J}/\mathcal{J}_N} (\hat{\mathbf{u}}_\alpha(t), \mathbf{v})_{\mathbb{L}_2}^2 / \alpha! = 0.$$

Thus, $(\bar{\mathbf{u}}(t), \mathbf{v})_{\mathbb{L}_2}$ is also a continuous process in $L_2(\Omega)$.

Now we shall prove that $\bar{\mathbf{u}}(t)$ is a global solution of (2.1). Let $Z$ be the set of real-valued sequences $z = (z_i^k)_{k \geq 1}$ such that only a finite number of $z_k$ are not zero. For $h(s, z) = z_i^k m_i(s) \ell_k$, write

$$p_t(z) = \exp \left\{ \int_0^t h(s, z) \cdot dW(s) - \tfrac{1}{2} \int_0^t |h(s, z)|_Y^2 \, ds \right\}.$$

It is readily seen (see, e.g., [24]) that

$$p_t(z) = \sum_{\alpha \in \mathcal{J}} \frac{z^\alpha}{\alpha!} \zeta_t^\alpha. \tag{3.15}$$

Obviously, (3.15) yields

$$\bar{\mathbf{u}}^z(t) := \mathbf{E}\bar{\mathbf{u}}(t) p_t(z) = \sum_{\alpha \in \mathcal{J}} \hat{\mathbf{u}}_\alpha(t) \frac{z^\alpha}{\alpha!}. \tag{3.16}$$

Thus, taking into account that $\{\hat{\mathbf{u}}_\alpha(t),\ \alpha \in \mathcal{J}\}$ verifies (3.6), it is easily checked that for every $t \in (0,T]$ and every test function $\varphi \in \mathcal{V}$,

$$\int_{\mathbf{R}^d} (\bar{\mathbf{u}}^z(t), \varphi)\, dx$$
$$= I_{\{|\alpha|=0\}} \int_{\mathbf{R}^d} (\mathbf{u}(0), \varphi)\, dx$$
$$+ \int_0^t \int_{\mathbf{R}^d} \Bigg\{ (a^{ij}(s)\, \partial_j \bar{\mathbf{u}}^z(s), \partial_i \varphi)$$
$$+ \Bigg( \varphi, -\sum_{\alpha \in \mathcal{J}} z^\alpha \sum_{\gamma, \beta \in U^\alpha} (\hat{u}^i_\gamma(s)\varphi, \partial_i \hat{\mathbf{u}}_\beta(s)) \Phi(\alpha, \beta, \gamma)$$
$$+ b^i(s) \partial_i \bar{\mathbf{u}}^z(s) + \tilde{L}_i(\bar{\mathbf{u}}^z(s), t) \mathbf{h}^i(s)$$
$$+ I_{\{|\alpha|=0\}} (\mathbf{f}(t) + \partial_j \mathbf{f}^j(t, \mathbf{u}(t))) \Bigg)$$
$$+ \Bigg( m_j(s) \alpha^k_j \varphi, \sigma^{ik}(t) \sum_{\alpha \in \mathcal{J}} \partial_i \hat{\mathbf{u}}_{\alpha(j,k)}(s) \frac{z^\alpha}{\alpha!}$$
$$+ \sum_{\alpha \in \mathcal{J}} I_{\{|\alpha|=1\}} \mathbf{g}(s) \frac{z^\alpha}{\alpha!} \Bigg) \Bigg\} dx\, ds,$$

**P**-a.a.

Owing to (3.15),

$$\mathbf{E} p_t(z) \int_0^t \int_{\mathbf{R}^d} (\varphi, \bar{u}^k(s)\, \partial_k \bar{\mathbf{u}}(s))\, dx\, ds$$
$$= \sum_{\alpha \in \mathcal{J}} \frac{z^\alpha}{\alpha!} \int_0^t \int_{\mathbf{R}^d} (\varphi, \mathbf{E}(\bar{u}^k(s)\, \partial_k \bar{\mathbf{u}}(s) \zeta_\alpha))\, dx\, ds.$$

Now, by mimicking the proof of (3.13), one can show

$$\mathbf{E} p_t(z) \int_0^t \int_{\mathbf{R}^d} (\varphi, \bar{u}^k(s)\, \partial_k \bar{\mathbf{u}}(s))\, dx\, ds$$
$$= \sum_{\alpha \in \mathcal{J}} z^\alpha \int_0^t \int_{\mathbf{R}^d} \sum_{\gamma, \beta \in U^\alpha} (\hat{u}^i_\gamma(s)\varphi, \partial_i \hat{\mathbf{u}}_\beta(s)) \Phi(\gamma, \beta, \alpha)\, dx\, ds.$$

Next, let us consider

$$F^i(t) := \sum_{\alpha \in \mathcal{J}} \frac{z^\alpha}{\alpha!} \int_0^t \int_{\mathbf{R}^d} \sum_{j,k} (\varphi, \sigma^{ik}(t)\, \partial_i \hat{\mathbf{u}}_{\alpha(j,k)}(s)) m_j(s) \alpha_j^k \, dx \, ds.$$

Denote $\psi^i(s) := (\varphi, \sum_i \sigma^i(t)\, \partial_i \bar{\mathbf{u}}(s))$ and write $\psi^{ik} = (\psi^i, \ell_k)_Y$, $\psi^{ik}_{\alpha(j,k)}(s) = \mathbf{E}\psi^{ik}(s)\zeta_{\alpha(j,k)}$ (no summation over $k$ is assumed here). It is readily checked that

$$F^i(t) = \sum_{\alpha \in \mathcal{J}} \frac{z^\alpha}{\alpha!} \int_0^t \int_{\mathbf{R}^d} \sum_{j,k} \psi^{ik}_{\alpha(j,k)}(s) m_j(s) \alpha_j^k \, dx \, ds$$

(3.17)
$$= \int_0^t \int_{\mathbf{R}^d} \sum_{j,k} z_j^k \sum_{\alpha \in \mathcal{J}} \frac{z^{\alpha(j,k)}}{\alpha(j,k)!} \psi^{ik}_{\alpha(j,k)}(s) m_j(s) \, dx \, ds$$

$$= \mathbf{E}\left[ p_t(z) \int_0^t \sum_{j,k} z_j^k m_j(s) \int_{\mathbf{R}^d} \psi^{ik}(s) \, dx \, ds \right].$$

On the other hand, by the Itô formula and (3.17),

$$\mathbf{E}\left[ p_t(z) \int_0^t \int_{\mathbf{R}^d} \psi^i(s) \, dx \, dW_s \right] = \mathbf{E}\left[ p_t(z) \int_0^t \sum_{j,k} z_j^k m_j(s) \int_{\mathbf{R}^d} \psi^{ik}(s) \, dx \, ds \right]$$

$$= F^i(t).$$

Similarly, one can prove

$$\int_0^t \int_{\mathbf{R}^d} \left( \varphi, \sum_{\alpha \in \mathcal{J}} I_{\{|\alpha|=1\}} \mathbf{g}(s) \frac{z^\alpha}{\alpha!} \right) dx \, ds$$

$$= E p_t(z) \int_0^t \int_{\mathbf{R}^d} (\varphi, g(s)) \, dx \, dW_s.$$

Thus, we have proved that for every $t \in (0, T]$ and every test function $\varphi \in \mathcal{V}$,

$$E p_t(z) \int_{\mathbf{R}^d} (\bar{\mathbf{u}}(t), \varphi) \, dx$$

(3.18)
$$= E p_t(z) \Bigg\{ I_{\{|\alpha|=0\}} \int_{\mathbf{R}^d} (\mathbf{u}(0), \varphi) \, dx$$

$$+ \int_0^t \int_{\mathbf{R}^d} [(a^{ij}(s)\, \partial_j \bar{\mathbf{u}}(s), \partial_i \varphi)$$

$$+ (\varphi, -\bar{u}^k(s)\, \partial_k \bar{\mathbf{u}}(s) + b^i(s)\, \partial_i \bar{\mathbf{u}}^z(s)$$

$$+ \tilde{L}_i(\bar{\mathbf{u}}^z(s), t)\mathbf{h}^i(s)$$

$$+ I_{\{|\alpha|=0\}}(\mathbf{f}(t) + \partial_j \mathbf{f}^j(t, \mathbf{u}(t))))] \, dx \, ds$$

$$+ \int_0^t \int_{\mathbf{R}^d} \left( \varphi, \sum_i \sigma^i(t)\, \partial_i \bar{\mathbf{u}}(s) + \mathbf{g}(s) \right) dx \, dW_s \Bigg\}$$

$P$-a.a.

Now, it is not difficult to show that the linear subspace generated by $\{p_t(z)\}_z$ is dense in $L_2(\Omega, \mathcal{F}_t, P)$. Thus, $\bar{\mathbf{u}}(t)$ is a strong solution of (2.1). $\square$

Now we can prove that the uniqueness of a solution of (2.1) is equivalent to the uniqueness of a solution of (3.6).

THEOREM 3.5. *A strong $\mathbb{L}_2$-weakly continuous solution of* (2.1) *with the property* (3.3) *is pathwise unique if and only if a solution of* (3.6) *is unique in the class of $\mathbb{L}_2$-weakly continuous functions $\{\hat{\mathbf{u}}_\alpha(t,x), \alpha \in \mathcal{J}\}$ such that* (3.5) *holds.*

PROOF. Under the assumptions of Theorem 3.2, any strong solution of (2.1) is given by

$$\mathbf{u}(t) = \sum_{\alpha \in \mathcal{J}} \frac{\hat{\mathbf{u}}_\alpha(t)}{\alpha!} \zeta_\alpha.$$

Therefore, the uniqueness of a solution of the propagator equation yields pathwise uniqueness of a strong solution of (2.1).

If $\{\hat{\mathbf{u}}^1_\alpha(t,x), \alpha \in \mathcal{J}\}$ and $\{\hat{\mathbf{u}}^2_\alpha(t,x), \alpha \in \mathcal{J}\}$ are two solutions of equation (3.6) such that inequality (3.5) holds, then, by Theorem 3.4, a strong weakly continuous solution of (2.1) with the property (3.3) is pathwise unique if and only if a solution of (3.6) is unique in the class of $\mathbb{L}_2$-valued weakly continuous functions $\{\hat{\mathbf{u}}_\alpha(t,x), \alpha \in \mathcal{J}\}$ such that (3.5) holds.

Let $\{\hat{\mathbf{u}}^1_\alpha(t,x), \alpha \in \mathcal{J}\}$ and $\{\hat{\mathbf{u}}^2_\alpha(t,x), \alpha \in \mathcal{J}\}$ be two solutions of (3.6) with the property (3.5). Then, by Theorem 3.4, $\mathbf{u}^1(t) = \sum_{\alpha \in \mathcal{J}} \hat{\mathbf{u}}^1_\alpha(t) \zeta_\alpha / \alpha!$ and $\mathbf{u}^2(t) = \sum_{\alpha \in \mathcal{J}} \hat{\mathbf{u}}^2_\alpha(t) \zeta_\alpha / \alpha!$ are strong solutions of (2.1) with the property (3.3). The uniqueness of the strong solution of (2.1) yields that for every $\alpha$,

$$\hat{\mathbf{u}}^1_\alpha(t) = \mathbf{E}\mathbf{u}^1(t)\zeta_\alpha = \mathbf{E}\mathbf{u}^2(t)\zeta_\alpha = \hat{\mathbf{u}}^2_\alpha(t). \qquad \square$$

One could give another proof of Theorem 2.2 using Theorems 3.4 and 3.5 (in this regard, see Remark 3.1). We expect that this approach would be useful in proving the existence of a strong global solution of the Navier–Stokes equation (2.1) in the three-dimensional case.

REMARK 3.1. If $u = (u^l)_{1 \leq l \leq d}$ is an $\mathbb{L}_2$-weakly continuous solution of (2.1) such that

$$\sup_{s \leq T} \mathbf{E}|\mathbf{u}(s)|_2^2 + \int_0^T \mathbf{E}|\nabla \mathbf{u}(s)|_2^2 \, ds < \infty,$$

then

$$\mathbf{E} \sup_{s \leq T} |\mathbf{u}(s)|_2^2 < \infty.$$

Indeed, let $\varphi \in C_0^\infty(\mathbf{R}^d), \varphi \geq 0, \int \varphi \, dx = 1, \varphi_\varepsilon(x) = \varepsilon^{-d}\varphi(x/\varepsilon)$. Applying the Itô formula for $|u_\varepsilon(t)|_2^2 = |u(t) * \varphi_\varepsilon|_2^2$, we have

$$|\mathbf{u}_\varepsilon(t)|_2^2 = |\mathbf{u}_\varepsilon(0)|_2^2$$
$$+ \int_0^t \int \{2[-a^{ij}(s)\,\partial_i u^l(s) - f^{l,j}(s, \mathbf{u}(s))]_\varepsilon \,\partial_j u_\varepsilon^l(s)$$
$$+ 2[b^{l,k}(s)\,\partial_k u^l(s) + f^l(s, \mathbf{u}(s))$$
$$+ \tilde{L}^k(s, \mathbf{u}) \cdot h^{l,k}(s)]_\varepsilon u_\varepsilon^l(s)$$
$$+ |(\sigma^j(s)\,\partial_j u^l(s) + g^l(s, \mathbf{u}(s)) - \tilde{L}^l(\mathbf{u}, s))_\varepsilon|_Y^2\} \, dx \, ds$$
$$+ 2\int_0^t \int [\sigma^j(s)\,\partial_j u^l(s) + g^l(s, \mathbf{u}(s))]_\varepsilon u_\varepsilon^l(s) \, dx \, dW_s,$$

where

$$[-a^{ij}(s)\,\partial_i u^l(s) - f^{l,j}(s, \mathbf{u}(s))]_\varepsilon$$
$$= [-a^{ij}(s)\,\partial_i u^l(s) - f^{l,j}(s, \mathbf{u}(s))] * \varphi_\varepsilon,$$
$$= [b^{l,k}(s)\,\partial_k u^l(s) + f^l(s, \mathbf{u}(s)) + \tilde{L}^k(s, \mathbf{u}) \cdot h^{l,k}(s)]_\varepsilon$$
$$= [b^{l,k}(s)\,\partial_k u^l(s) + f^l(s, \mathbf{u}(s)) + \tilde{L}^k(s, \mathbf{u}) \cdot h^{l,k}(s)] * \varphi_\varepsilon,$$
$$(\sigma^j(s)\,\partial_j u^l(s) + g^l(s, \mathbf{u}(s)) - \tilde{L}^l(\mathbf{u}, s))_\varepsilon$$
$$= (\sigma^j(s)\,\partial_j u^l(s) + g^l(s, \mathbf{u}(s)) - \tilde{L}^l(\mathbf{u}, s)) * \varphi_\varepsilon,$$
$$= [\sigma^j(s)\,\partial_j u^l(s) + g^l(s, \mathbf{u}(s))]_\varepsilon$$
$$= [\sigma^j(s)\,\partial_j u^l(s) + g^l(s, \mathbf{u}(s))] * \varphi_\varepsilon.$$

Since, for each $\eta > 0$,

$$\mathbf{E}\sup_{t \leq T}\left|\int_0^t \int [\sigma^j(s)\,\partial_j u^l(s) + g^l(s, \mathbf{u}(s))]_\varepsilon u_\varepsilon^l(s) \, dx \, dW_s\right|$$
$$\leq \eta\mathbf{E}\sup_{t \leq T}|\mathbf{u}_\varepsilon(\mathbf{t})|_2^2 + C_\eta \mathbf{E}\int_0^T |\nabla \mathbf{u}(s)|_2^2 + |\mathbf{g}(s, \mathbf{u}(s))|_2^2 \, ds,$$

we derive easily that

$$\mathbf{E}\sup_{s \leq T} |\mathbf{u}(s)|_2^2 < \infty.$$

## APPENDIX

For the reader's convenience, we summarize below some useful facts regarding Helmholtz's decomposition of vector fields (see, e.g., [13, 23]).

The Riesz transform will be used for the definition of the projections. For $f \in L_2(\mathbf{R}^d, Y)$, set

$$R_j(f)(x) = \lim_{\varepsilon \to 0} c_* \int_{|y| \geq \varepsilon} \frac{y_j}{|y|^{d+1}} f(x-y)\, dy, \qquad j = 1, \ldots, d,$$

with $c_* = G(\frac{n+1}{2})/\pi^{(n+1)/2}$ ($G$ is the gamma function). $R_j$ is called a Riesz transform. According to [31], Chapter III, formula (8), page 58,

$$(R_j f)\hat{}\,(x) = -i \frac{\xi_j}{|\xi|} \cdot \hat{f},$$

where

$$\hat{f}(\xi) = \mathcal{F}(f) = (2\pi)^{-d/2} \int e^{-i(\xi, x)} f(x)\, dx.$$

Given a function $f \in L_p(\mathbf{R}^d, Y)$, we define a vector Riesz transform $Rf = (R_1 f, \ldots, R_d f)$.

For $\mathbf{v} \in \mathbb{L}_2(Y)$, set

$$\mathcal{G}(\mathbf{v}) = -RR_j v^j, \qquad \mathcal{S}(\mathbf{v}) = \mathbf{v} - \mathcal{G}(\mathbf{v}).$$

Then (see Lemma 2.7 in [23]), $\mathbb{L}_2(Y)$ is a direct sum

$$\mathbb{L}_2(Y) = \mathcal{G}(\mathbb{L}_2(Y)) \oplus \mathcal{S}(\mathbb{L}_2(Y)),$$

$$\mathcal{S}(\mathbb{L}_2(Y)) = \{\mathbf{g} \in \mathbb{L}_2(Y) : \operatorname{div} \mathbf{g} = 0\},$$

and $\mathcal{G}(\mathbb{L}_2(Y))$ is a Hilbert subspace orthogonal to $\mathcal{S}(\mathbb{L}_2(Y))$.

The functions $\mathcal{G}(\mathbf{v})$ and $\mathcal{S}(\mathbf{v})$ are usually referred to as the potential and the divergence-free projections, respectively, of the vector field $\mathbf{v}$.

The following statement holds.

LEMMA A.1 (see Lemmas 2.11 and 2.12 in [23]). *$\mathcal{G}, \mathcal{S}$ can be extended continuously to all $\mathbb{H}_2^s(Y)$, $s \in (-\infty, \infty)$: there is a constant $C$ so that for all $\mathbf{v} \in \mathbb{H}_2^s(Y)$,*

$$\|\mathcal{G}(\mathbf{v})\|_{s,2} \leq C \|\mathbf{v}\|_{s,2}, \qquad \|\mathcal{S}(\mathbf{v})\|_{s,2} \leq C \|\mathbf{v}\|_{s,2}.$$

*Moreover, the space $\mathbb{H}_2^s(Y)$ can be decomposed into the direct sum*

$$\mathbb{H}_2^s(Y) = \mathcal{G}(\mathbb{H}_2^s(Y)) \oplus \mathcal{S}(\mathbb{H}_2^s(Y)),$$

*and, if $\mathbf{f} \in \mathcal{G}(\mathbb{H}_2^s(Y))$, $\mathbf{g} \in \mathcal{S}(\mathbb{H}_2^{-s}(Y))$, then*

(A.1) $$\langle \mathbf{f}, \mathbf{g} \rangle_{\mathbb{H}_2^s(Y), \mathbb{H}_2^{-s}(Y)} = \langle \mathbf{f}, \mathbf{g} \rangle_s = 0.$$

*Also,*

(A.2) $$\mathcal{S}(\mathbb{H}_2^s(Y)) = \{\mathbf{v} \in \mathbb{H}_2^s(Y) : \operatorname{div} \mathbf{v} = \mathbf{0}\}.$$

**Acknowledgment.** We are indebted to the referee for insightful suggestions.

Department of Mathematics
  and Center of Applied Mathematical Sciences
University of Southern California
Los Angeles, California 90089
USA
e-mail: rozovski@math.usc.edu